\documentclass[11pt]{article}
\usepackage{amssymb,latexsym,amsmath,amsbsy,amsthm,amsxtra,amsgen}
\newfont{\msbm}{msbm10 at 11pt}

\newcommand {\Z} {\mbox{\msbm Z}}

\newcommand {\N} {\mbox{\msbm N}}

\newfont{\msbmsm}{msbm10 at 8pt}

\newcommand {\Zsm} {\mbox{\msbmsm Z}}

\newfont{\msbmmed}{msbm10 at 10pt}

\newcommand {\Zmed} {\mbox{\msbmmed Z}}

\newtheorem{Theo}{Theorem}[section]
\newtheorem{Lemma}[Theo]{Lemma}
\newtheorem{Cor}[Theo]{Corollary}
\newtheorem{Prop}[Theo]{Proposition}

\newtheorem{Rmk}[Theo]{Remark}

\begin{document}
\title{The loop-erased random walk and the uniform spanning tree on the four-dimensional discrete torus}

\author{by Jason Schweinsberg\thanks{Supported in part by NSF Grant DMS-0504882}}
\maketitle

\begin{abstract}
Let $x$ and $y$ be points chosen uniformly at random from $\Zmed_n^4$, the four-dimensional discrete torus with side length $n$.  We show that the length of the loop-erased random walk from $x$ to $y$ is of order $n^2 (\log n)^{1/6}$, resolving a conjecture of Benjamini and Kozma.
We also show that the scaling limit of the uniform spanning tree on $\Zmed_n^4$ is the Brownian continuum random tree of Aldous.  Our proofs use the techniques developed by Peres and Revelle, who studied the scaling limits of the uniform spanning tree on a large class of finite graphs that includes the $d$-dimensional discrete torus for $d \geq 5$, in combination with results of Lawler concerning intersections of four-dimensional random walks.
\end{abstract}

\footnote{{\it AMS 2000 subject classifications}.  Primary 60K35;
Secondary 60G50, 60D05}

\footnote{{\it Key words and phrases}.  Loop-erased random walk, Uniform spanning tree, Continuum random tree.}

\section{Introduction}

Given a finite connected graph $G$, a spanning tree of $G$ is a subgraph of $G$ which contains every vertex of $G$ and has no cycles.  Every connected graph $G$ has at least one spanning tree, and a spanning tree picked uniformly at random from the set of all possible spanning trees is called a uniform spanning tree.  In this paper, we study the scaling limit of the uniform spanning tree on the graph $\Z^4_n = \{0, 1, \dots, n-1\}^4$, the four-dimensional discrete torus of side length $n$, as $n \rightarrow \infty$.

Previously, Aldous \cite{CRT1} studied the scaling limit of the uniform spanning tree on the complete graph on $m$ vertices as $m \rightarrow \infty$.  He showed that the limit is an object called the continuum random tree (CRT), which we now describe using the line-breaking construction in section 4 of \cite{CRT3}.  Consider an inhomogeneous Poisson process on $[0, \infty)$ with rate $r(t) = t$, and denote the points of this Poisson process by $0 < t_1 < t_2 < \dots$.  To construct the CRT, begin with a single line segment of length $t_1$, and label the endpoints $z_1$ and $z_2$.  Then choose a random point on this segment, attach to it another line segment of length $t_2 - t_1$, and label $z_3$ the endpoint of this new segment that is not on the line connecting $z_1$ and $z_2$.  To continue the process inductively, once $k$ line segments have been added to the tree, add an additional line segment of length $t_{k+1} - t_k$ to a uniformly chosen point on the existing tree, and label the endpoint of the new line segment that was not on the previous tree $z_{k+2}$.  We view each line segment as being orthogonal to all previous line segments, and the limiting object that we get as $k \rightarrow \infty$, called the CRT, is a random metric space.  The points $z_1, z_2, \dots$ can be viewed as points chosen uniformly at random from a ``mass measure" on the CRT.  If $d(z_i, z_j)$ denotes the distance between the points $z_i$ and $z_j$, then we denote by $\mu_k$ the joint distribution of the distances $(d(z_i, z_j))_{1 \leq i < j \leq k}$.

To state more precisely Aldous' result for the uniform spanning tree on the complete graph, let $K_m$ denote the complete graph on $m$ vertices, and let $\tilde{\cal T}$ be a uniform spanning tree on $K_m$.  Fix a positive integer $k$, and let $y_1, \dots, y_k$ be vertices chosen independently and uniformly at random from $K_m$.  Let $\tilde{d}(y_i, y_j)$ be the distance between $y_i$ and $y_j$, that is, the number of vertices along the path in $\tilde{\cal T}$ from $y_i$ to $y_j$.  Then, as $m \rightarrow \infty$,
\begin{equation}
\bigg( \frac{\tilde{d}(y_i, y_j)}{\sqrt{m}} \bigg)_{1 \leq i < j \leq k} \rightarrow_d \mu_k.
\label{aldCG}
\end{equation}

Recently, Peres and Revelle \cite{perrev} proved a conjecture of Pitman by showing that the CRT also arises as the limit of the uniform spanning tree on a large class of finite graphs that includes the hypercubes $\Z^2_n$, expander graphs, and the $d$-dimensional torus $\Z^d_n$ when $d \geq 5$.  Parts of their proof, however, do not work in dimension four.  Our main result, stated below, is that the CRT is also the scaling limit of the uniform spanning tree on $\Z^4_n$.

\begin{Theo}
Fix a positive integer $k$, and let $x_1, \dots, x_k$ be points chosen independently and uniformly at random from the four-dimensional discrete torus $\Z^4_n$.  Let ${\cal T}$ be a uniform spanning tree on $\Z^4_n$, and let $d(x_i, x_j)$ be the number of vertices along the path in ${\cal T}$ from $x_i$ to $x_j$.  There exists a sequence of constants $(\gamma_n)_{n=1}^{\infty}$, bounded away from zero and infinity, such that the joint distribution of the distances $$\frac{d(x_i, x_j)}{\gamma_n n^2 (\log n)^{1/6}},$$ with $1 \leq i < j \leq k,$ converges to $\mu_k$ as $n \rightarrow \infty$.
\label{mainth}
\end{Theo}

Note that Theorem \ref{mainth} implies that the typical distance in the uniform spanning tree between two randomly chosen points in $\Z^4_n$ is of order $n^2 (\log n)^{1/6}$, as compared with $n^{d/2}$ on $\Z^d_n$ with $d \geq 5$.  The exponent of $1/6$ in the logarithmic correction was conjectured by Benjamini and Kozma \cite{benkoz}. 
In view of the connections between the uniform spanning tree and the loop-erased random walk, to be reviewed in the next subsection, and the fact that Brownian motion is the scaling limit of loop-erased random walk on $\Z^4$, it is natural that the CRT is the scaling limit of the uniform spanning tree on the torus in dimension four.  We do not expect to get the CRT in the limit in dimensions below four.  See Lawler, Schramm, and Werner \cite{lsw} for a discussion of the uniform spanning tree in two dimensions, and Kozma \cite{kozma} for recent progress on the scaling limit of the loop-erased random walk in three dimensions.  See also \cite{schw} for results on the dynamics of the loop-erased random walk on the torus in dimensions four and higher which are based on developments in the present paper.

\subsection{Loop-erased random walks and Wilson's algorithm}

We recall here the definition of a loop-erased random walk.  A path in a graph $G$ is a sequence of vertices $u_0, u_1, \dots, u_j$ such that $u_{i-1}$ and $u_i$ are connected by an edge for $i = 1, \dots, j$.  Given a path $\lambda$ in $G$, let $|\lambda|$ denote the number of vertices in the path $\lambda$.  Also, let $LE(\lambda)$ denote the loop-erasure of $\lambda$, obtained by erasing loops from the path $\lambda$ in chronological order.  More formally, if $\lambda = (u_0, u_1, \dots, u_j)$, then $LE(\lambda)$ is the path $(v_0, \dots, v_k)$ obtained inductively as follows.  First set $v_0 = u_0$.  Suppose $v_0, \dots, v_m$ have been defined for some $m \geq 0$.
If $v_m = u_j$, then $k = m$ and $v_m$ is the last vertex in the path $LE(\lambda)$.  Otherwise, define $v_{m+1} = u_{r+1}$, where $r = \max\{i: u_i = v_m\}$.

A random walk on a graph $G$ is a Markov chain $X = (X_t)_{t=0}^{\infty}$ taking its values in the set of vertices of $G$ such that $P(X_{t+1} = w|X_t = v) = d(v)^{-1} {\bf 1}_{\{v \sim w\}}$, where we write $d(v)$ for the degree of the vertex $v$ and $v \sim w$ whenever $v$ and $w$ are connected by an edge.  That is, at each step the random walk moves to a neighboring vertex of the graph, chosen uniformly at random.  We denote the random walk segment $(X_t)_{t=a}^b$ by $X[a,b]$.  Given a vertex $x$ and a subset of vertices $V$, a loop-erased random walk from $x$ to $V$ is the path $LE(X[0,T])$, where $X_0 = x$ and $T = \min\{t: X_t \in V\}$.  At times we will consider random walks that stay in their current position with probability $1/2$ and move to a randomly chosen neighboring vertex with probability $1/2$.  However, it is clear from the definition that repeating vertices along a path does not affect the loop erasure of the path.

It is well-known that there are strong connections between uniform spanning trees and loop-erased random walks.  Pemantle \cite{pemantle} showed that if $x$ and $y$ are vertices in a finite graph $G$, then the loop-erased random walk from $x$ to $y$ has the same distribution as the path from $x$ to $y$ in a uniform spanning tree of $G$.  Pemantle's result, in combination with Theorem \ref{mainth}, immediately gives the following corollary.

\begin{Cor}
Let $x$ and $y$ be points chosen independently and uniformly at random from $\Z^4_n$.  Let $X = (X_t)_{t=0}^{\infty}$ be a random walk on $\Z^4_n$ with $X_0 = x$, and let $T = \min\{t: X_t = y\}$.  Then there is a sequence of constants $(\gamma_n)_{n=1}^{\infty}$, bounded away from zero and infinity, such that for all $x \geq 0$,
$$\lim_{n \rightarrow \infty} P \big( |LE(X[0,T])| > \gamma_n n^2 (\log n)^{1/6} x \big) = e^{-x^2/2}.$$
\label{Rayleighcor}
\end{Cor}

\noindent The limiting distribution in Corollary \ref{Rayleighcor} is known as the Rayleigh distribution.  It arises because the probability that the points $z_1$ and $z_2$ in the CRT are a distance at least $x$ apart is the probability that there are no points in $[0, x]$ in an inhomogeneous Poisson process of rate $r(t) = t$, which is $e^{-x^2/2}$.

Wilson \cite{wilson} established an even stronger connection between uniform spanning trees and loop-erased random walks by discovering an algorithm for constructing uniform spanning trees using loop-erased random walks.  The algorithm proceeds as follows.  First pick some vertex $x_0$ in the graph $G$, and let $T_0 = \{x_0\}$.  Then given a tree $T_i$ for some $i \geq 0$, choose a vertex $x_{i+1}$ in $G$ and define the tree $T_{i+1}$ by adjoining to $T_i$ the loop-erased random walk from $x_{i+1}$ to $T_i$.  Because loops are erased, each $T_i$ is a tree, and if we continue the process until $T_i$ contains every vertex of $G$, the tree $T_i$ is a uniform spanning tree.  We emphasize that the algorithm works for any choice of the vertices $x_0, x_1, \dots$, with the choice of $x_{i+1}$ even being permitted to depend on the tree $T_i$.  Our proof of Theorem \ref{mainth}, which we outline in the next subsection, will use Wilson's algorithm to compare uniform spanning trees on $\Z^4_n$ and on the complete graph $K_m$.

\subsection{Outline of the proof}

In this subsection, we outline the main ideas in the proof of Theorem \ref{mainth}.  The structure of the proof is borrowed from the proof of Peres and Revelle \cite{perrev} for $d \geq 5$.  The idea is to couple a uniform spanning tree on $\Z^4_n$ with a uniform spanning tree on the complete graph, so that Theorem \ref{mainth} will follow from this coupling and Aldous' result that the CRT arises as the scaling limit of the uniform spanning tree on the complete graph.  However, some of the steps in the proof involve applying results about intersections of four-dimensional random walks which are more intricate than the results needed to push through the proof when $d \geq 5$. 

Let $x_1, \dots, x_k$ be points chosen independently and uniformly at random from $\Z^4_n$, and let $y_1, \dots, y_k$ be points chosen independently and uniformly at random from the complete graph $K_m$.  We can construct a uniform spanning tree on $K_m$ by using Wilson's algorithm, starting the first $k$ random walks from the points $y_1, \dots, y_k$.  Note that for the purposes of studying the joint distribution of the distances between pairs of these points, we can stop the algorithm after $k$ steps because once $y_i$ and $y_j$ are in the tree, the distance between them will not change when other branches are added to the tree.  Likewise, we can construct a uniform spanning tree on $\Z^4_n$ by using Wilson's algorithm and starting the first $k$ walks from $x_1, \dots, x_k$.

To couple these two spanning trees, we break the random walks on $\Z^4_n$ into segments, with each segment corresponding to a single step of a random walk on $K_m$.  Then, for the coupling, we collapse the vertices in the spanning tree on $\Z^4_n$ that come from the same segment of a random walk to a single vertex.  There is some flexibility in the choice of the length of these segments, but to obtain the best bound for the probability that the coupling succeeds, we choose the length of the segments to be $$r = \lfloor n^2 (\log n)^{9/22} \rfloor.$$  Within a segment of length $r$, the random walk will make many short loops that get erased during the construction.  Long loops, in which one segment intersects the loop-erasure of a previous segment, can cause entire segments to get erased.

On $K_m$, assuming there is a self-loop at every vertex, two steps of a random walk coincide with probability $1/m$.  Therefore, for a long loop to be equally likely to form in both constructions, enabling the coupling to succeed, the probability that a random walk segment on $\Z^4_n$ of length $r$ intersects the loop-erasure of another such segment must be approximately $1/m$.  Two random walks on $\Z^4_n$ of length $r$ intersect with probability of order $r^2/(n^4 \log n) \approx (\log n)^{-2/11}$ (see Propositions \ref{mainintprop} and \ref{lowerintprop} below), and from results of Lyons, Peres, and Schramm \cite{lps03}, it will follow that the probability that a random walk of length $r$ intersects a loop-erased walk of length $r$ is the same order of magnitude.  Therefore, we need to choose $m$ to be of order $(\log n)^{2/11}$.  More precisely, following \cite{perrev}, we define the capacity of a set $U \subset \Z^4_n$ to be
\begin{equation}
\mbox{Cap}_r(U) = P(X_t \in U \mbox{ for some }t \leq r),
\label{capdef}
\end{equation}
where $(X_t)_{t=0}^{\infty}$ is a random walk on $\Z^4_n$ started at a uniformly chosen point.
We will show in Proposition \ref{maincapprop} that the capacity of a loop-erased random walk of length $r$ is tightly concentrated around its mean of $a_n (\log n)^{-2/11}$, for some constant $a_n$.  Therefore, if we choose $m = \lfloor a_n^{-1} (\log n)^{2/11} \rfloor$, then the probability that the next segment of length $r$ intersects the loop-erasure of a given previous segments is always very close to $1/m$, enabling the coupling with the process on the complete graph to hold with high probability.  Because $r$ is much longer than the mixing time of a random walk on $\Z^4_n$, the fact that the starting point of each segment is not chosen independently at random does not have a large effect on these estimates.

When the coupling of spanning trees holds, the number of vertices along the path in the spanning tree on $K_m$ from $y_i$ to $y_j$ is the same as the number of segments of length $r$ along the path in the spanning tree on $\Z^4_n$ from $x_i$ to $x_j$.  To couple the distances, it is necessary to estimate the length of the paths of length $r$ after loop-erasure, and to show that the distribution of this length is highly concentrated around its mean.  It follows from a result of Lawler \cite{law95} that the fraction of the points on each segment that are retained after loop-erasure is of order $(\log n)^{-1/3}$, so the length of the segments after loop-erasure is of order $r/(\log n)^{1/3}$.  We will show in Proposition \ref{lengthprop}, using arguments in \cite{lawler}, that this length is concentrated around its mean.

There is one further complication that must be addressed.  To carry out Wilson's algorithm on $\Z^4_n$, it is necessary to start the tree with an initial vertex.  Then we start a random walk from $x_1$ and run the walk until it hits that vertex.  However, the time before a random walk on $\Z^4_n$ hits a given vertex is of order $n^4$, much longer than the length of a loop-erased random walk between two points.  To avoid having to run a random walk for a time of order $n^4$, we again use a technique of Peres and Revelle \cite{perrev} by adding a root vertex to each graph.  Given $\alpha > 0$, let $K_{m, \alpha}$ be a weighted graph constructed by starting with $K_m$, giving every edge of $K_m$ (including the self-loop at every vertex) weight $1$, and then adding a root vertex $\rho$ that is connected to every other vertex in the graph by an edge of weight $m/(\alpha \sqrt{m} - 1)$.  This means a weighted random walk on $K_{m, \alpha}$ that on each step moves to a neighboring vertex with probabilities proportional to the edge weights is the same as the simple random walk on $K_m$, except that it jumps to $\rho$ after a number of steps which is geometrically distributed with mean $\alpha \sqrt{m}$.  Likewise, let $G_{n, \beta}$ be obtained by adding to $\Z_n^4$ a root vertex $\rho$, which is connected to every other vertex in the graph by an edge of weight $8/(\beta n^2 (\log n)^{1/2} - 1)$.  The weighted random walk on $\Z^4_n$ is the same as the simple random walk on $\Z^4_n$, until it jumps to $\rho$ after a number of steps which is geometrically distributed with mean $\beta n^2 (\log n)^{1/2}$.  Given $\beta$, we choose 
$\alpha$ so that 
\begin{equation}
\frac{1}{\alpha \sqrt{m}} = 1 - \bigg(1 - \frac{1}{\beta n^2 (\log n)^{1/2}} \bigg)^r.
\label{alphadef}
\end{equation}
This choice ensures that the probability that a single step of a random walk on $K_{m, \alpha}$ visits the root is the same as the probability that an $r$-step random walk on $G_{n, \beta}$ visits the root.  It is then possible to couple uniform spanning trees on $K_{m, \alpha}$ and $G_{n, \beta}$ by starting the tree with the root vertex, so that the initial random walks started from $y_1$ and $x_1$ respectively are run until they hit the root.  Once we have uniform spanning trees on $K_{m, \alpha}$ and $G_{n, \beta}$, we can remove the edges leading to the root to get spanning forests on $K_m$ and $\Z^4_n$.  As $\beta \rightarrow \infty$, the probability that $x_1, \dots, x_n$ are in the same tree component of the spanning forest tends to one, so we can use results on stochastic domination of spanning forests by spanning trees, as in \cite{perrev}, to couple the spanning trees on $K_m$ and $\Z^4_n$.

We conclude this discussion by observing that, with the above picture in mind, there is a simple explanation for why the exponent $1/6$ arises in the logarithmic correction.  If we run two random walks started from $x_1$ and $x_2$ for time $L \geq n \log n$, then they intersect with probability of order $\min\{1, L^2/(n^4 \log n)\}$, which is of order one when $L = n^2 (\log n)^{1/2}$.  Therefore, the path in the uniform spanning tree between $x_1$ and $x_2$ comes from random walk paths whose lengths are of order $n^2 (\log n)^{1/2}$.  The fraction of points along these paths that survive loop-erasure is of order $(\log n)^{-1/3}$, so the lengths of the remaining paths are of order $n^2 (\log n)^{1/6}$.

The rest of this paper is organized as follows.  In section 2, we prove some results about random walks on the four-dimensional torus that will be needed to study loop-erased random walks.  In section 3, we study loop-erased walks on the torus, focusing on obtaining tight bounds on the length and capacity of a loop-erased segment.  Then in section 4, we use these results to couple a uniform spanning tree on $\Z^4_n$ with a uniform spanning tree on $K_m$ and complete the proof of Theorem \ref{mainth}.

\section{Random walks on the four-dimensional torus}

In this section, we establish some facts about random walks on $\Z_n^4$ that we will need to study the loop-erased random walk and the uniform spanning tree.  Most of the key ideas come from the book of Lawler \cite{lawler}, which contains numerous results for random walks on $\Z^4$.  The work in this section primarily involves establishing the analogous results for random walks on the torus. 

We denote the point $(0, 0, 0, 0)$ by $0$, and if $x = (x_1, \dots, x_4) \in \Z_n^4$, we denote the Euclidean norm by $|x| = (x_1^2 + \dots + x_4^2)^{1/2}$.  It will be convenient to work with aperiodic random walks.  Therefore, when we say that $X = (X_t)_{t=0}^{\infty}$ is a random walk on $\Z_n^4$, we will assume that at each step, the random walk stays in its current position with probability $1/2$.  That is, for all $t$, we have $P(X_{t+1} = X_t|X_t) = 1/2$ and, if $x$ is one of the $8$ points in $\Z^4$ such that $|x| = 1$, then $P(X_{t+1} = X_t + x|X_t) = 1/16$.  To apply results of Lawler \cite{lawler} about random walks on $\Z^4$, it will be necessary to consider also random walks on $\Z^4$ that never stay in their current position.  Therefore, if $Z = (Z_t)_{t=0}^{\infty}$ is a random walk on $\Z^4$, we will specify that $Z$ is a simple random walk if $P(Z_{t+1} = Z_t|Z_t) = 0$ and a lazy random walk if $P(Z_{t+1} = Z_t|Z_t) = 1/2$.  In both cases, when the random walk does not stay in its current position, it moves to one of the $8$ neighboring points with equal probability.  All of our random walks will be in discrete time.

Throughout the paper, $C, C', C_1, C_2, \dots$ will denote positive constants that do not depend on $n$ but whose values may change from line to line.

\subsection{Bounds on transition probabilities}
\label{retsubsec}

In this subsection, we will establish some bounds on transition probabilities for random walks on $\Z^4_n$.  Let $X = (X_t)_{t=0}^{\infty}$ be a random walk on $\Z^4_n$.  Let $Y = (Y_t)_{t=0}^{\infty}$ be a simple random walk on $\Z^4$, and let $Z = (Z_t)_{t=0}^{\infty}$ be a lazy random walk on $\Z^4$.  For all $x, y \in \Z_n^4$, let $p_{t,n}(x, y) = P(X_t = y|X_0 = x)$.  For all $x, y \in \Z^4$, let $q'_t(x, y) = P(Y_t = y|Y_0 = x)$ and let $q_t(x, y) = P(Z_t = y|Z_0 = x)$.  The following result for the simple random walk, which is a version of the Local Central Limit Theorem, is the four-dimensional case of Theorem 1.2.1 in \cite{lawler}.

\begin{Lemma}
There exists a constant $C$ such that for all $t \in \N$ and $x = (x_1, \dots, x_4) \in \Z^4$ such that $t + x_1 + \dots + x_4$ is even, we have $$\bigg| q'_t(0, x) - \frac{8}{\pi^2 t^2} e^{-2 |x|^2/t} \bigg| \leq C \min \bigg\{ \frac{1}{t^3}, \frac{1}{t^2 |x|^2} \bigg\}.$$
\label{LawlerLCLT}
\end{Lemma}

The following large deviations result, known as Azuma's Inequality, will be useful for deducing properties of the lazy random walk from properties of the simple random walk.  If $(M_t)_{t=0}^{\infty}$ is a martingale with $M_0 = 0$ and $|M_{t+1} - M_t| \leq c$ for all $t$, then
\begin{equation}
P(M_t \geq x) \leq e^{-x^2/2c^2t}.
\label{azuma}
\end{equation}
See section 2.4 of \cite{jlr} for a proof.  Azuma's Inequality and Lemma \ref{LawlerLCLT} lead to the following bounds for transition probabilities of the lazy random walk.

\begin{Lemma}
There exists a constant $C_1$ such that
\begin{equation}
q_t(0, 0) \leq \frac{C_1}{t^2}
\label{lazyreturn}
\end{equation}
for all $t \in \N$.  There also exists a constant $C_2$ such that for all $t \in \N$ and all $x \in \Z^4$, we have
\begin{equation}
q_t(0, x) \leq \frac{C_2}{|x|^4}.
\label{lazyLD}
\end{equation}
\end{Lemma}

\begin{proof}
It follows from the $x = 0$ case of Lemma \ref{LawlerLCLT} that there exists a constant $C$ such that $q'_t(0,0) \leq C t^{-2}$ for all $t \in \N$.  To obtain the bound for the lazy random walk, we condition on the number $N_t$ of steps $s \leq t$ such that $Z_s \neq Z_{s-1}$.  If $M_t = t/2 - N_t$, then $(M_t)_{t=0}^{\infty}$ is a martingale with $|M_{t+1} - M_t| \leq 1/2$ for all $t$.  Therefore $P(N_t < t/4) \leq e^{-t/8}$ by (\ref{azuma}).  It follows that $q_t(0,0) \leq C (t/4)^{-2} + e^{-t/8}$, which implies (\ref{lazyreturn}).

To prove (\ref{lazyLD}), note that if $a > 0$, the function $t \mapsto t^{-2} e^{-2a/t}$, defined for $t \in (0, \infty)$, is maximized when $t = a$.  Therefore, $t^{-2} e^{-2|x|^2/t} \leq e^{-2} |x|^{-4}$ for all $t \in \N$ and $x \in \Z^4$.  Since $q_t'(0, x) = 0$ when $t < |x|$ and $1/(t^2 |x|^2) \leq |x|^{-4}$ when $t \geq |x|$, we get $q_t'(0, x) \leq C |x|^{-4}$.  Because this bound does not depend on $t$, we get the same bound for $q_t(0, x)$ by conditioning on $N_t$.
\end{proof}

The next step is to use these bounds on the transition probabilities for the random walk on $\Z^4$ to obtain bounds on the transition probabilities for the random walk on the torus $\Z^4_n$.  Note that the random walks $X$ and $Z$ can be coupled so that $X_t \equiv Z_t$ (mod $n$) for all $t$.  

\begin{Prop}
Suppose $x = (x_1, \dots, x_4)$ is a point in $\Z^4$ such that $|x_i| \leq n/2$ for $i =1, \dots, 4$.  Let $x' = (x_1', \dots, x_4') \in \Z_n^4$ be defined so that $x_i' \equiv x_i$ (mod $n$) for $i = 1, \dots, 4$.  Then for any constant $C'$, there is a constant $C$ such that for all $t \leq C' n^2$, we have
\begin{equation}
q_t(x, 0) \leq p_{t,n}(x', 0) \leq q_t(x, 0) + \frac{C}{n^4}.
\label{returnbound}
\end{equation}
\label{torustrans}
\end{Prop}

\begin{proof}
Let $V = (n \Z)^4$.  It is clear from the coupling between $X$ and $Z$ that
\begin{equation}
p_{t,n}(x', 0) = \sum_{y \in V} q_t(x, y).
\label{sumqt}
\end{equation}
The first inequality in (\ref{returnbound}) follows immediately.  To obtain the second inequality, we divide $\Z^4$ into boxes.  Given $y = (y_1, \dots, y_4) \in \Z^4$, let $$B(y) = \{(z_1, \dots, z_4): |y_i - z_i| \leq n/2 \mbox{ for } i = 1, \dots, 4\}.$$
Denote by $\partial B(y)$ the set of all points in $B(y)$ having a neighbor that is not in $B(y)$.

Assume that $Z_0 = x$, so $q_t(x, y) = P(Z_t = y)$.  Suppose $y \in V \setminus \{0\}$.  If $Z_t = y$, then $Z_s \in \partial B(y)$ for some $s \leq t$.  Because $|y - z| \geq (n-1)/2$ for all $z \in \partial B(y)$, equation (\ref{lazyLD}) implies there is a constant $C$ such that for all $s \leq t$ and all $z \in \partial B(y)$ we have $P(Z_t = y|Z_s = z) \leq C n^{-4}$.  By applying the strong Markov property at the stopping time $\min\{t: Z_t \in B(y)\}$, it follows that $q_t(x, y) \leq C n^{-4} P(Z_s \in B(y) \mbox{ for some }s \leq t)$.  Therefore, it remains to show that there is a constant $C$ such that $$\sum_{y \in V \setminus \{0\}} P(Z_t \in B(y) \mbox{ for some }s \leq t) \leq C.$$

For $k \geq 1$, the cardinality of $$V_k = \bigg\{y = (y_1, \dots, y_4) \in V: \frac{kn}{2} \leq \max_{1 \leq i \leq 4} |x_i - y_i| < \frac{(k+1)n}{2} \bigg\}$$ is at most $(k+1)^4$.  For $y \in V_k$, there is an $i \in \{1, \dots, 4\}$ such that $|x_i - z_i| \geq (k-1)n/2$ for all $z = (z_1, \dots, z_4) \in B(y)$.  Let $Z_{i,t}$ be the $i$th coordinate of $Z_t$, and note that $(Z_{i,t})_{t=0}^{\infty}$ is a martingale.  By applying the Reflection Principle followed by (\ref{azuma}) to $(Z_{i,t})_{t=0}^{\infty}$, we get
$$P(Z_t \in B(y) \mbox{ for some }s \leq t) \leq 2 P(|Z_{i,t}| \geq (k-1)n/2) \leq 4 e^{-(k-1)^2 n^2/8t},$$
and the right-hand side is at most $4  e^{-(k-1)^2/8C'}$ whenever $t \leq C'n^2$.  It follows that
$$\sum_{y \in V \setminus \{0\}} P(Z_t \in B(y) \mbox{ for some }s \leq t) \leq \sum_{k=1}^{\infty} 4(k+1)^4 e^{-(k-1)^2/8C'} \leq C,$$ which implies the proposition.
\end{proof}

\begin{Cor}
For any constant $C'$, there exists a constant $C$ such that $p_{t,n}(0, 0) \leq C/t^2$ for all $t$ such that $1 \leq t \leq C' n^2$.
\label{returntorus}
\end{Cor}

\begin{proof}
Since $1/n^4 \leq (C'/t)^2$, the result follows immediately from (\ref{lazyreturn}) and the upper bound in Proposition \ref{torustrans}.
\end{proof}

Proposition \ref{torustrans} and Corollary \ref{returntorus} apply only when $t \leq C' n^2$.  To bound transition probabilities when $t > C' n^2$, we will need facts about the mixing time of the random walk $X$, which are discussed in the next subsection.

\subsection{Mixing time}

The random walk $X$ is an irreducible, aperiodic Markov chain whose stationary distribution is uniform on $\Z_n^4$.  Therefore (see, for example, chapter 5 of \cite{durrett}), the distribution of $X_t$ converges as $t \rightarrow \infty$ to the uniform distribution on $\Z^4_n$, which we denote by $\pi$.  To measure the speed of convergence, we define the uniform mixing time
$$\tau_n = \min \bigg\{t: \sup_{x \in \Zsm^4_n} \bigg| p_{t,n}(0, x) - \frac{1}{n^4} \bigg| \leq \frac{1}{2n^4} \bigg\}.$$  For any two probability measures $\mu$ and $\nu$ on $\Z_n^4$, define $$d(\mu, \nu) = \sup_{x \in \Zsm_n^4} n^{4} |\mu(x) - \nu(x)|,$$ so $\tau_n = \min\{t: d(p_{t,n}(0, \: \cdot), \pi) \leq 1/2\}$.  A standard subadditivity result on convergence for Markov chains (see, for example, the $p = \infty$ case of Proposition 2.2 in \cite{saloffcoste}) implies that if $d(p_{t,n}(0, \: \cdot), \pi) \leq \beta$ and $s \geq t$, then $d(p_{s,n}(0, \: \cdot), \pi) \leq \beta^{\lfloor s/t \rfloor}$.  In particular, for all $t \geq \tau_n$, we have
\begin{equation}
\sup_{x \in \Zsm^4_n} \bigg| p_{t,n}(0, x) - \frac{1}{n^4} \bigg| \leq n^{-4} 2^{-\lfloor t/\tau_n \rfloor}.
\label{tightmixbound}
\end{equation}
It is well-known that there exist constants $0 < C_1 < C_2 < \infty$ such that
\begin{equation}
C_1 n^2 \leq \tau_n \leq C_2 n^2 \hspace{.2in} \mbox{ for all }n.
\label{mixeq}
\end{equation}
See, for example, the calculations in chapter 5 of \cite{aldfill}.  The notion of mixing time will be useful in what follows primarily because of the following lemma.

\begin{Lemma}
\label{mixlem}
Let $f$ be a nonnegative function of two paths in $\Z_n^4$ of lengths $L+1$ and $M+1$.  Suppose $X = (X_t)_{t=0}^{\infty}$ and $Y = (Y_t)_{t=0}^{\infty}$ are independent random walks on $\Z_n^4$, and suppose $Y_0$ has the uniform distribution on $\Z_n^4$.  If $d \geq \tau_n$, then
$$E[f(X[0,L], X[L+d, L+M+d])] \leq \frac{3}{2} E[f(X[0,L], Y[0,M])].$$
\end{Lemma}

\begin{proof}
Let $x = (x_0, \dots, x_L)$ be a path in $\Z_n^4$ of length $L+1$, and let $y = (y_0, \dots, y_M)$ be a path in $\Z^4_n$ of length $M+1$.  By the definition of $\tau_n$,
$$P(X_{L+d} = y_0|X[0,L] = x) \leq \max_{z \in \Zsm_n^4} p_{d,n}(x,z) = \max_{z \in \Zsm_n^4} p_{d,n}(0,z) \leq \frac{3}{2n^4} = \frac{3}{2} P(Y_0 = y_0|X[0,L] = x).$$  Therefore, by the Markov property,
\begin{align}
P(X[0,L] = x&, X[L+d, L+M+d] = y) \nonumber \\
&= P(X[0,L] = x) P(X[L+d, L+M+d] = y|X[0,L] = x) \nonumber \\
&\leq \frac{3}{2} P(X[0,L] = x) P(Y[0,M] = y|X[0,L] = x) \nonumber \\
&= \frac{3}{2} P(X[0,L] = x, Y[0,M] = y), \nonumber
\end{align}
and the result follows easily.
\end{proof}

\subsection{Intersections of two walks}

Our goal in this subsection is to prove upper and lower bounds for the probability that two random walks on the torus intersect.  Our first lemma concerns the expected number of intersections of two random walks started at the origin.

\begin{Lemma}
Let $(X_t)_{t=0}^{\infty}$ and $(Y_t)_{t=0}^{\infty}$ be two independent random walks on $\Z_n^4$ started at the origin.  Let $$R = \sum_{s=0}^L \sum_{t=0}^L {\bf 1}_{\{X_s = Y_t\}}$$ be the number of intersections.  If $L \leq C' n^2 (\log n)^{1/2}$ for some constant $C'$, then there exists a constant $C$ such that $E[R] \leq C \log n$.
\label{originlem}
\end{Lemma}

\begin{proof}
Note that $P(X_s = Y_t) = P(X_{s + t} = 0)$.  It follows from (\ref{mixeq}) and Corollary \ref{returntorus} that there is a constant $C$ such that $P(X_t = 0) \leq C/t^2$ if $1 \leq t \leq \tau_n$.  If $t \geq \tau_n$, then $P(X_t = 0) \leq 3/2n^4$.  Therefore,
\begin{align}
E[R] &= \sum_{s=0}^L \sum_{t=0}^L P(X_s = Y_t) \leq \sum_{t = 0}^{\tau_n} (t + 1) P(X_t = 0) + \frac{3 (L+1)^2}{2 n^4} \nonumber \\
&\leq 1 + \sum_{t=1}^{\tau_n} \frac{C (t+1)}{t^2} + \frac{3 (L+1)^2}{2 n^4} \leq C \bigg( \log \tau_n + \frac{L^2}{n^4} \bigg), \nonumber
\end{align}
which gives the result.
\end{proof}

For all $x \in \Z_n^4$, denote by $G_L(x)$ the expected number of times that a random walk of length $L$ started at $0$ visits $x$.  That is, if $(X_t)_{t=0}^{\infty}$ is a random walk on $\Z_n^4$ started at $0$, then $$G_L(x) = \sum_{t=0}^{L} P(X_t = x).$$  Likewise, for $x \in \Z^4$, let $G_L^*(x)$ be the expected number of times that a simple random walk on $\Z^4$ of length $L$ started at $0$ visits $x$.  Let $G^*(x)$ be the expected number of times that a simple random walk on $\Z^4$ of infinite length started at $0$ visits $x$, so $G^*$ is the standard Green's function for random walks on $\Z^4$.

\begin{Lemma}
Let $X = (X_t)_{t=0}^{\infty}$ be a random walk on $\Z_n^4$ started at the origin.  Let $$D_n = \sum_{t=0}^n G_n(X_t).$$  There exist constants $C_1$ and $C_2$ such that $P(D_n \leq C_1 \log n) \leq C_2/\log n$ for all $n$.
\label{Dnlem}
\end{Lemma}

\begin{proof}
Let $Z = (Z_t)_{t=0}^{\infty}$ be a lazy random walk on $\Z^4$, coupled with $X$ so that $X_t \equiv Z_t$ (mod $n$) for all $t$.  Let $Z' = (Z'_t)_{t=0}^{\infty}$ be another lazy random walk on $\Z^4$ started at $0$, independent of $Z$, and let $X' = (X'_t)_{t=0}^{\infty}$ be the corresponding random walk on $\Z_n^4$ such that $X'_t \equiv Z'_t$ (mod $n$) for all $t$.  Since $X_t = X'_s$ whenever $Z_t = Z'_s$, we have
\begin{equation}
D_n = E \bigg[ \sum_{t=0}^n \sum_{s=0}^n {\bf 1}_{\{X_t = X'_s\}} \bigg| X \bigg] \geq E \bigg[ \sum_{t=0}^n \sum_{s=0}^n {\bf 1}_{\{Z_t = Z'_s\}} \bigg| Z \bigg].
\label{Dndef}
\end{equation}
Let $T_0 = T'_0 = 0$, and for $k \geq 1$, let $T_k = \min\{t > T_{k-1}: Z_t \neq Z_{T_{k-1}}\}$ and $T'_k = \min\{t > T'_{k-1}: Z'_t \neq Z'_{T_{k-1}}\}$.  Then let $Y_t = Z_{T_t}$ and $Y'_t = Z'_{T'_t}$ for all $t \geq 0$, so $Y = (Y_t)_{t=0}^{\infty}$ and $Y' = (Y'_t)_{t=0}^{\infty}$ are simple random walks on $\Z_n^4$ started at $0$.

Let $$D^* = \sum_{t=0}^{\lfloor n/4 \rfloor} G^*(Y_t).$$  Lawler showed that there exist constants $C_3$ and $C_4$ such that
\begin{equation}
P(D^* \leq C_3 \log n) \leq \frac{C_4}{\log n}
\label{lawc3c4}
\end{equation}
(see p. 101 of \cite{lawler}, or use Proposition 3.4.1 of \cite{lawler} in combination with Chebyshev's Inequality).  Note that $q'_s(0, 0) \leq C/s^2$ for all $s \geq 1$ by Lemma \ref{LawlerLCLT}.  Therefore, $$E \bigg[ \sum_{t=0}^{\lfloor n/4 \rfloor} (G^*(Y_t) - G_{\lfloor n/4 \rfloor}^*(Y_t)) \bigg] = \sum_{t=0}^{\lfloor n/4 \rfloor} \sum_{s=\lfloor n/4 \rfloor+1}^{\infty} P(Y_t = Y'_s) \leq (\lfloor n/4 \rfloor +1) \sum_{s=\lfloor n/4 \rfloor+1}^{\infty} q_s'(0, 0) \leq C.$$  It follows that if 
\begin{equation}
D_n^* = \sum_{t=0}^{\lfloor n/4 \rfloor} G_{\lfloor n/4 \rfloor}^*(Y_t) = E \bigg[ \sum_{t=0}^{\lfloor n/4 \rfloor} \sum_{s=0}^{\lfloor n/4 \rfloor} {\bf 1}_{\{Y_s = Y'_s\}} \bigg| Y \bigg],
\label{Dnstardef}
\end{equation}
then $P(D^* - D_n^* \geq (C_3 \log n)/ 2) \leq C/(\log n)$ by Markov's Inequality.  Combining this observation with (\ref{lawc3c4}), we get
\begin{equation}
P \bigg( D_n^* \leq \frac{C_3 \log n}{2} \bigg) \leq \frac{C}{\log n}.
\label{Dnstarlog}
\end{equation}

Azuma's Inequality (\ref{azuma}) gives
$P(T_{\lfloor n/4 \rfloor} > n) \leq e^{-n/8}$ and $P(T'_{\lfloor n/4 \rfloor} > n) \leq e^{-n/8}$.  Therefore, by (\ref{Dndef}) and (\ref{Dnstardef}), we have $D_n \geq D_n^* - (1 + n/4)^2 e^{-n/8}$ on the event that $T_{\lfloor n/4 \rfloor} \leq n$.  The lemma now follows from (\ref{Dnstarlog}).
\end{proof}

Our main result concerning intersection probabilities of two random walks on the torus is the following.  We prove it by using essentially the same proof used by Lawler to bound the intersection probability of two random walks in $\Z^4$ (see the proof of the upper bound of Theorem 3.3.2 in \cite{lawler}).

\begin{Prop}
Let $(X_t)_{t=0}^{\infty}$ and $(Y_t)_{t=0}^{\infty}$ be independent random walks on $\Z_n^4$.  Assume that $Y_0$ has the uniform distribution on $\Z^4_n$.  Also, assume there is a constant $C'$ such that $K \geq C' n \log n$ and $L \geq C' n \log n$.  Then there is a constant $C$ such that $$P(X_s = Y_t \mbox{ for some }s \leq K \mbox{ and }t \leq L) \leq \frac{CKL}{n^4 \log n}.$$
\label{mainintprop}
\end{Prop}

\begin{proof}
Set $\kappa = \min\{t \geq 0: X_t = Y_s \mbox{ for some }s \leq L\}$, and let $\sigma = \min\{s \geq 0: X_{\kappa} = Y_s\}$.  We need to bound $P(\kappa \leq K)$.  Let $$J = \bigg\{ s: \sum_{t=0}^n G_n(Y_{t+s} - Y_s) \leq C_1 \log n \bigg\},$$ where $C_1$ is the constant from Lemma \ref{Dnlem}.  For all $t = 0, 1, \dots, K$ and $s = 0, 1, \dots, L$, we have $P(X_t = Y_s) = 1/n^4$ because $Y_s$ has the stationary distribution $\pi$.  Also, the event that $s \in J$ is independent of the event that $X_t = Y_s$.  By Lemma \ref{Dnlem}, we have $P(s \in J) \leq C_2/(\log n)$ whenever $s \leq L-n$.  Putting these results together, we get
\begin{align}
P(\kappa \leq K, \sigma \in J) &= \sum_{s=0}^L P(\kappa \leq K, \sigma = s, s \in J) \nonumber \\
&\leq \sum_{t=0}^K \sum_{s=0}^{L-n} P(X_t = Y_s, s \in J) + \sum_{t=0}^K \sum_{s=L-n+1}^L P(X_t = Y_s) \nonumber \\
&\leq \sum_{t=0}^K \sum_{s=0}^{L-n} \frac{C_2}{n^4 \log n} + \sum_{t=0}^K \sum_{s=L-n+1}^L \frac{1}{n^4} \nonumber \\
&\leq \frac{C_2 (K + 1)L}{n^4 \log n} + \frac{K+1}{n^3}.
\label{intbound1}
\end{align}
Also,
\begin{equation}
P(K-n < \kappa \leq K) \leq \sum_{t=K-n+1}^K \sum_{s=0}^L P(X_t = Y_s) = \frac{n(L + 1)}{n^4}.
\label{intbound2}
\end{equation}
Let $$R = \sum_{t=0}^K \sum_{s=0}^L {\bf 1}_{\{X_t = Y_s\}}.$$  By the definition of $J$, we can apply the strong Markov property to $(X_t)_{t=0}^K$ while holding the path $(Y_t)_{t=0}^{\infty}$ fixed to get
$$E[R|\kappa \leq K-n, \sigma \notin J] \geq C_1 \log n.$$  It follows that
\begin{equation}
P(\kappa \leq K-n, \sigma \notin J) \leq \frac{E[R]}{C_1 \log n} = \frac{1}{C_1 \log n} \sum_{t=0}^K \sum_{s=0}^L P(X_t = Y_s) = \frac{(K+1)(L+1)}{C_1 n^4 \log n}.
\label{intbound3}
\end{equation}
The result now follows from (\ref{intbound1}), (\ref{intbound2}), and (\ref{intbound3}).
\end{proof}

The next result gives a lower bound for the intersection probability when the two walks have the same length.  

\begin{Prop}
Let $X = (X_t)_{t=0}^{\infty}$ and $Y = (Y_t)_{t=0}^{\infty}$ be independent random walks on $\Z^4_n$.  Assume that $Y_0$ has the uniform distribution on $\Z^4_n$.  If there is a constant $C'$ such that $L \leq C' n^2 (\log n)^{1/2}$, then there is a constant $C$ such that 
$$P(X_s = Y_t \mbox{ for some }s \leq L \mbox{ and }t \leq L) \geq \frac{C L^2}{n^4 \log n}.$$
\label{lowerintprop}
\end{Prop}

In the proof of this proposition, as well as later in the paper, we will use the following well-known result (see, for example, (3.5) on p. 95 of \cite{lawler}), which gives a lower bound on the probability that a nonnegative random variable is strictly positive.

\begin{Lemma}
Let $Z$ be a nonnegative random variable.  Then $P(Z > 0) \geq (E[Z])^2/E[Z^2]$.
\label{2moment}
\end{Lemma}

\begin{proof}[Proof of Proposition \ref{lowerintprop}]
Let $$R = \sum_{s=0}^L \sum_{t=0}^L {\bf 1}_{\{X_s = Y_t\}}.$$  Since $Y_0$ has the uniform distribution, $P(X_s = Y_t) = 1/n^4$ for all $s$ and $t$, so $$E[R] = \frac{(L+1)^2}{n^4}.$$  Also, if $s \leq u$, then by the Markov property and time reversibility of $Y$, we have $$P(X_u = Y_v|X_s = Y_t) = P(X_{u-s} = Y_{|v-t|}|X_0 = Y_0).$$  Thus,
\begin{align}
E[R^2] &= E \bigg[ \sum_{s=0}^L \sum_{t=0}^L \sum_{u=0}^L \sum_{v=0}^L {\bf 1}_{\{X_s = Y_t\}} {\bf 1}_{\{X_u = Y_v\}} \bigg] \nonumber \\
&\leq 2 \sum_{s=0}^L \sum_{t=0}^L \sum_{u=s}^L \sum_{v=0}^L P(X_s = Y_t) P(X_u = X_v|X_s = Y_t) \nonumber \\
&= 2 \sum_{s=0}^L \sum_{t=0}^L \sum_{u=s}^L \sum_{v=0}^L \frac{1}{n^4} P(X_{u-s} = Y_{|v-t|}|X_0 = Y_0) \nonumber \\
&= 4 \sum_{s=0}^L \sum_{t=0}^L \sum_{q=0}^L \sum_{r=0}^L \frac{1}{n^4} P(X_q = Y_r|X_0 = Y_0) \nonumber \\
&= \frac{4(L+1)^2}{n^4} \sum_{q=0}^L \sum_{r=0}^L P(X_q = Y_r|X_0 = Y_0) \leq \frac{C L^2 \log n}{n^4}. \nonumber
\end{align}
where the last inequality follows from Lemma \ref{originlem}.  By Lemma \ref{2moment}, $$P(R > 0) \geq \frac{(E[R])^2}{E[R^2]} \geq \frac{C L^2}{n^4 \log n},$$ and the result follows.
\end{proof}

\subsection{Self-intersections}

We now use Proposition \ref{mainintprop} to bound the probability that a random walk on $\Z^4_n$ intersects its previous path.  These bounds will be important for the study of the loop-erased random walk on $\Z^4_n$.  Proposition \ref{medloops} below shows that random walks whose length is much shorter than $n^2 (\log n)^{1/2}$ are unlikely to make loops whose length is much larger than $n^2$.

\begin{Prop}
Let $(X_t)_{t=0}^{\infty}$ be a random walk on $\Z^4_n$.  There exists a constant $C$ such that $$P(X_t = X_s \mbox{ for some }s \leq L \mbox{ and }t \leq L \mbox{ such that }t \geq s + 2 \tau_n) \leq \frac{C L^2}{n^4 \log n}.$$
\label{medloops}
\end{Prop}

\begin{proof}
We may assume $L \geq 2 \tau_n$.  Let $m = \lceil (L+1)/\tau_n \rceil$, and for $i = 1, \dots, m$, let $B_i = \{(i-1)\tau_n, \dots, i \tau_n - 1\}$.
If $s + 2 \tau_n \leq t \leq L$, then there exist $B_i$ and $B_j$ with $j \geq i + 2$ such that $s \in B_i$ and $t \in B_j$.  If $j \geq i + 2$, then by Lemma \ref{mixlem} and Proposition \ref{mainintprop}, $$P(X_t = X_s \mbox{ for some }s \in B_i \mbox{ and }t \in B_j) \leq \frac{3 C \tau_n^2}{2 n^4 \log n},$$ where $C$ is the constant from Proposition \ref{mainintprop}.  Since the number of pairs of segments $B_i$ and $B_j$ is at most $m^2$, it follows that $$P(X_t = X_s \mbox{ for some }s \mbox{ and }t \mbox{ such that }t \geq s + 2 \tau_n) \leq
\bigg( \frac{L+1}{\tau_n} + 1 \bigg)^2 \frac{3 C \tau_n^2}{2 n^4 \log n},$$ which gives the result.
\end{proof}

Our next result uses the coupling introduced in subsection \ref{retsubsec} between a random walk on $\Z_n^4$ and a random walk on $\Z^4$.  We show that a random walk on $\Z_n^4$ whose length is much shorter than $n^2 (\log n)^{1/2}$ is unlikely to have loops that do not also appear in the corresponding walk on $\Z^4$.  This is closely related to Proposition \ref{medloops} because typically short loops are present in random walks on $\Z^4$ while long loops occur only on the torus.

\begin{Prop}
Let $(X_t)_{t=0}^{\infty}$ be a random walk on $\Z^4_n$.  Let $(Z_t)_{t=0}^{\infty}$ be a lazy random walk on $\Z^4$.  Assume the random walks are coupled so that $X_t \equiv Z_t$ (mod $n$) for all $t$.
There is a constant $C$ such that
$$P(X_s = X_t \mbox{ and }Z_s \neq Z_t \mbox{ for some }s \leq L, t \leq L) \leq \frac{CL^2}{n^4 \log n}.$$
\label{torusintprop}
\end{Prop}

\begin{proof}
The proof is similar to the proof of Proposition \ref{mainintprop}. 
Define $\kappa = \min\{t \geq 0: X_s = X_t \mbox{ and }Z_s \neq Z_t \mbox{ for some }s < t\}$, so $\kappa$ is the first time that there is a loop in the walk on $\Z^4_n$ that does not correspond to a loop in the walk on $\Z^4$.  Let $\sigma = \min\{s \geq 0: X_s = X_{\kappa} \mbox{ and }Z_s \neq Z_{\kappa}\}$.  We need to bound $P(\kappa \leq L)$.  Note that if $X_s = X_t$ and $t - s < n$, then $Z_s = Z_t$.  Therefore, if $\kappa \leq L$, then $\sigma \leq \kappa - n$.  Define $$J = \bigg\{s: \sum_{t=0}^{\lfloor n/3 \rfloor} G_{\lfloor n/3 \rfloor}(X_{t+s} - X_s) \leq C_1 \log (\lfloor n/3 \rfloor) \bigg\},$$ where $C_1$ is the constant from Lemma \ref{Dnlem}.

Proposition \ref{torustrans} implies that that there is a constant $C$ such that $p_{t, n}(0,0) - q_t(0,0) \leq C/n^4$ whenever $1 \leq t < \tau_n$, while (\ref{tightmixbound}) and the definition of $\tau_n$ imply that if $t \geq \tau_n$, then we have $p_{t, n}(0,0) - q_t(0,0) \leq p_{t, n}(0,0) \leq 3/2n^4$.  Therefore, there is a constant $C$ such that if $s < t$, then
\begin{equation}
P(X_s = X_t \mbox{ and }Z_s \neq Z_t) = p_{t-s, n}(0, 0) - q_{t-s}(0,0) \leq \frac{C}{n^4}.
\label{transdiff}
\end{equation}
It follows that
\begin{equation}
P(L - \lfloor n/3 \rfloor \leq \kappa \leq L) \leq \sum_{s=0}^L \sum_{t = \max\{s+1, L - \lfloor n/3 \rfloor\}}^L P(X_s = X_t, Z_s \neq Z_t) \leq \frac{CL}{n^3}.
\label{kappaprob2}
\end{equation}

Let $$R = \sum_{s=0}^L \sum_{t=s+1}^L {\bf 1}_{\{X_s = X_t\}} {\bf 1}_{\{Z_s \neq Z_t\}}.$$
By (\ref{transdiff}), there is a constant $C$ such that $E[R] \leq CL^2/n^4$.
Note that if $X_s = X_t$ and $Z_s \neq Z_t$, then if $X_{s+u} = X_{t+v}$ for some $u$ and $v$ such that $u + v < n$, we have $Z_{s+u} \neq Z_{t+v}$.  It follows by applying the strong Markov property at time $\kappa$ that $$E[R| \kappa \leq L - \lfloor n/3 \rfloor, \sigma \notin J] \geq C_1 \log (\lfloor n/3 \rfloor).$$  Therefore,
\begin{equation}
P(\kappa \leq L - \lfloor n/3 \rfloor, \sigma \notin J) \leq \frac{E[R]}{C_1 \log (\lfloor n/3 \rfloor)} \leq \frac{C L^2}{n^4 \log n}.
\label{kappaprob3}
\end{equation}

By Lemma \ref{Dnlem}, there is a constant $C$ such that $P(s \in J) \leq C/(\log n)$ for all $s \leq L - n/3$.  Therefore,
\begin{align}
P(\kappa \leq L, \sigma \in J) &= \sum_{s=0}^{L-n} \sum_{t=s+n}^L P(\sigma = s, \kappa = t, s \in J) \nonumber \\
&\leq \sum_{s=0}^{L-n} \sum_{t=s+n}^L P(X_s = X_t, Z_s \neq Z_t, s \in J) \nonumber \\
&\leq \frac{C}{\log n} \sum_{s=0}^{L-n} \sum_{t=s+n}^L P(X_s = X_t, Z_s \neq Z_t|s \in J).
\label{kappaprob}
\end{align}

The next step is to use Proposition \ref{torustrans} to bound $P(X_s = X_t, Z_s \neq Z_t|s \in J)$.
The key idea is that the event $\{s \in J\}$ depends only on the random variables $Z_s, Z_{s+1}, \dots, Z_{s + \lfloor n/3 \rfloor}$, and so conditioning on this event does not have a large effect on the probability that $X_s = X_t$ and $Z_s \neq Z_t$.  Let $B = \{x = (x_1, \dots, x_4) \in \Z^4: |x_i| \leq \lfloor n/3 \rfloor \mbox{ for }i = 1, \dots, 4\}$.  Given $x \in \Z^4$, let $x' \in \Z_n^4$ be defined so that $x' \equiv x$ (mod $n$).  Therefore, using the Markov property at the time $s + \lfloor n/3 \rfloor$, we get that for $t \geq s + n$,
\begin{align}
P(X_s = X_t, Z_s \neq Z_t|s \in J) &\leq \max_{x \in \Zsm^4} \max_{y \in x + B} P(X_s = X_t, Z_s \neq Z_t|s \in J, Z_s = x, Z_{s + \lfloor n/3 \rfloor} = y) \nonumber \\
&= \max_{x \in \Zsm^4} \max_{y \in x + B} P(X_s = X_t, Z_s \neq Z_t|Z_s = x, Z_{s + \lfloor n/3 \rfloor} = y)
\nonumber \\
&= \max_{z \in B} \big( p_{t - s - \lfloor n/3 \rfloor, n}(z', 0) - q_{t-s-\lfloor n/3 \rfloor}(z,0) \big).
\label{condloopprob}
\end{align}
By Proposition \ref{torustrans}, there is a constant $C$ such that $p_{t - s - \lfloor n/3 \rfloor, n}(z', 0) - q_{t-s-\lfloor n/3 \rfloor}(z,0) \leq C/n^4$ whenever $z \in B$ and $t - s - \lfloor n/3 \rfloor < \tau_n$.  If $t - s - \lfloor n/3 \rfloor \geq \tau_n$ then $$p_{t - s - \lfloor n/3 \rfloor, n}(z', 0) - q_{t-s-\lfloor n/3 \rfloor}(z,0) \leq p_{t - s - \lfloor n/3 \rfloor, n}(z', 0) \leq 3/2n^4.$$  Combining these bounds with (\ref{kappaprob}) and (\ref{condloopprob}), we get
\begin{equation}
P(\kappa \leq L, \sigma \in J) \leq \frac{CL^2}{n^4 \log n}.
\label{kappaprob1}
\end{equation}
The proposition follows from (\ref{kappaprob2}), (\ref{kappaprob3}), and (\ref{kappaprob1}) when $L \geq n \log n$.  The result follows from Azuma's Inequality (\ref{azuma}) when $L < n \log n$ because in this case it takes longer than time $L$ for the random walk to move a distance $n$.
\end{proof}

\subsection{Local cutpoints}

Let $(X_t)_{t=0}^{\infty}$ be a random walk on $\Z_n^4$.  Using terminology similar to that of Peres and Revelle \cite{perrev}, we call the time $u$ a $k$-local cutpoint if $$\{X_{u - k}, \dots, X_{u-1} \} \cap \{X_{u+1}, \dots, X_{u + k}\} = \emptyset.$$  In other words, the portion of the random walk between times $u - k$ and $u-1$ does not intersect the portion between times $u+1$ and $u + k$.  We call a $2 \tau_n$-local cutpoint simply a local cutpoint.  Our goal in this subsection is to prove that for random walks on $\Z_n^4$, with high probability every sufficiently long segment contains at least one local cutpoint.  More precisely, we have the following proposition.

\begin{Prop}
\label{lcprop}
Let $C' > 0$ be a constant.  Let $(X_t)_{t=0}^{L}$ be a random walk on $\Z^4_n$, and suppose $L \leq C' n^2 (\log n)^{1/2}$.  Then for all $\theta > 0$, there is a constant $C$ such that with probability at least $1 - C/(\log n)$, every interval $[t, t + n^2 (\log n)^{\theta}]$ with $0 \leq t \leq L - n^2 (\log n)^{\theta}$ contains a local cutpoint.
\end{Prop}

Before proving Proposition \ref{lcprop}, we prove two lemmas.  The first gives an estimate on the probability that a given point is a local cutpoint.

\begin{Lemma}
Let $X = (X_t)_{t=0}^{\infty}$ and $X' = (X'_t)_{t=0}^{\infty}$ be independent random walks on $\Z^4_n$ starting at $0$.  For all $m$, let $f_n(m)$ be the probability that the two walks do not intersect through time $m$, so $$f_n(m) = P(X_s \neq X'_t \mbox{ whenever }1 \leq s \leq m, 1 \leq t \leq m).$$  Then $$\lim_{n \rightarrow \infty} \frac{f_n(2 \tau_n)}{f_n(n)} = \frac{1}{\sqrt{2}}.$$
\label{sqrt2lem}
\end{Lemma}

\begin{proof}
Let $Z = (Z_t)_{t=0}^{\infty}$ and $Z' = (Z'_t)_{t=0}^{\infty}$ be two lazy random walks on $\Z^4$ starting at $0$.  For all $m$, let $$g(m) = P(Z_s \neq Z'_t \mbox{ whenever }1 \leq s \leq m, 1 \leq t \leq m).$$  Lawler \cite{law92} showed that there is a constant $C_0$ such that
\begin{equation}
g(m) \sim \frac{C_0}{(\log m)^{1/2}}
\label{gmLawler}
\end{equation}
where $\sim$ means that the ratio of the two sides tends to $1$ as $m \rightarrow \infty$.  Lawler's estimate was stated for simple random walks, rather than lazy random walks.  However, (\ref{azuma}) implies that the number of steps on which a lazy random walk of length $m$ does not stay in the same place is between $m/4$ and $m$ with probability at least $1 - e^{-m/8}$, so it is easy to see that Lawler's estimate also holds for the lazy random walk.

We now claim that
\begin{equation} 
\lim_{n \rightarrow \infty} \frac{f_n(n)}{g(n)} = \lim_{n \rightarrow \infty} \frac{f_n(2\tau_n)}{g(2\tau_n)} = 1.
\label{fngeq}
\end{equation}
We may assume that the random walks are coupled so that $X_t \equiv Z_t$ (mod $n$) and $X_t' \equiv Z_t'$ (mod $n$) for all $t$.  Fix $m \leq 2 \tau_n$.  Define a random walk $(Y_t)_{t=0}^{2m}$ on $\Z^4$ by concatenating the walks $Z$ and $Z'$.  More precisely, define $Y_t = Z'_{m-t} - Z'_m$ for $0 \leq t \leq m$ and $Y_t = Z_{t - m} - Z'_m$ for $m+1 \leq t \leq 2m$.  If $Z_s = Z'_t$ for some $1 \leq s \leq m$ and $1 \leq t \leq m$, then $Y_{m-t} = Y_{m+s}$.  Now define a random walk $(W_t)_{t=0}^{2m}$ on $\Z_n^4$ such that $W_t \equiv V_t$ (mod $n$) for $0 \leq t \leq 2m$.  If $X_s = X'_t$ but $Z_s \neq Z'_t$ for some $1 \leq s \leq m$ and $1 \leq t \leq m$, then $Y_{m-t} \neq Y_{m+s}$ but $Y_{m-t} \equiv Y_{m+s}$ (mod $n$), so $W_{m-t} = W_{m+s}$.  Since $m \leq 2 \tau_n$, it follows from (\ref{mixeq}) and Proposition \ref{torusintprop} that $$P(W_s = W_t \mbox{ and }Y_s \neq Y_t \mbox{ for some }s < t \leq 2m) \leq \frac{C \tau_n^2}{n^4 \log n} \leq \frac{C}{\log n}.$$
Therefore, $0 \leq g(m) - f_n(m) \leq C/(\log n)$ for all $m \leq 2 \tau_n$.  In view of (\ref{gmLawler}), equation (\ref{fngeq}) follows.

Using (\ref{gmLawler}), (\ref{fngeq}), and (\ref{mixeq}), we get
$$\lim_{n \rightarrow \infty} \frac{f_n(2 \tau_n)}{f_n(n)} = \lim_{n \rightarrow \infty} \frac{g(2\tau_n)}{g(n)} = \lim_{n \rightarrow \infty} \frac{(\log 2 \tau_n)^{-1/2}}{(\log n)^{-1/2}} = \frac{1}{\sqrt{2}},$$
as claimed.
\end{proof}

\begin{Lemma}
Let $(X_t)_{t=0}^{\infty}$ be a random walk on $\Z_n^4$.  There exists a constant $b > 0$ such that if $L \leq 5 \tau_n$, the probability that $X[0,L]$ has a local cutpoint in $[2 \tau_n, L - 2 \tau_n]$ is at least $b$.
\label{newlclem}
\end{Lemma}

\begin{proof}
Choose points $x_1, \dots, x_k \in [2 \tau_n, L - 2 \tau_n]$ such that $x_{i+1} - x_i > 2n$ for $i = 1, \dots, k-1$.  Since $L \geq 5 \tau_n$, this can be done so that $k \geq C n$ for some positive constant $C$.  Let $A_i$ be the event that $x_i$ is an $n$-local cutpoint.  Since $x_{i+1} - x_i > 2n$, the events $A_1, \dots, A_k$ are independent, and each has probability $f_n(n)$.  Let $Y$ be the number of the points $x_1, \dots, x_k$ that are $n$-local cutpoints.  Then $Y$ has the binomial distribution with parameters $k$ and $f_n(n)$.  Let $Z$ be the number of the points $x_1, \dots, x_k$ that are local cutpoints.  Then $Z \leq Y$.  Since $$P(x \mbox{ is a local cutpoint}|x \mbox{ is an }n\mbox{-local cutpoint}) = \frac{f_n(2 \tau_n)}{f_n(n)},$$ 
we have, for $0 < a < 1/\sqrt{2}$, $$E[Z] = \frac{f_n(2 \tau_n)E[Y]}{f_n(n)} > a E[Y]$$ for sufficiently large $n$ by Lemma \ref{sqrt2lem}.  Using Lemma \ref{2moment}, we have, for sufficiently large $n$,
\begin{equation}
P(Z > 0) \geq \frac{(E[Z])^2}{E[Z^2]} \geq \frac{a^2 (E[Y])^2}{E[Y^2]} =
\frac{a^2 (k f_n(n))^2}{k f_n(n)(1 - f_n(n)) + (k f_n(n))^2}.
\label{PZ2mom}
\end{equation}
Note that $k f_n(n) \rightarrow \infty$ as $n \rightarrow \infty$ by (\ref{gmLawler}), (\ref{fngeq}), and the fact that $k \geq C n$.  Thus, the right-hand side of (\ref{PZ2mom}) approaches $a^2 > 0$ as $n \rightarrow \infty$, and the lemma follows.
\end{proof}

\begin{proof}[Proof of Proposition \ref{lcprop}]
Divide $[0, L]$ into disjoint subintervals of length $\frac{1}{2} n^2 (\log n)^{\theta}$, throwing away any remainder at the end of the interval.  If each of these subintervals has a local cutpoint, then so does every interval of length $n^2 (\log n)^{\theta}$.  Therefore, it suffices to show that with probability at least $1 - C/(\log n)$, each subinterval has a local cutpoint.  Note that the number of subintervals is at most $2C' (\log n)^{1/2 - \theta}$.

Now split an interval $I$ of length $\frac{1}{2} n^2 (\log n)^{\theta}$ into disjoint subintervals $J_1, \dots, J_m$ of length $5 \tau_n$, again throwing away any remainder.  By (\ref{mixeq}), $m \geq C_1 (\log n)^{\theta}$ for some constant $C_1 > 0$.  For $j = 1, \dots, m$, let $A_j$ be the event that there is a local cutpoint in the segment of length $\tau_n$ in the middle of $J_j$.  Since $A_j$ depends only on the random walk during times in $J_j$, the events $A_1, \dots, A_m$ are independent.  By Lemma \ref{newlclem}, $P(A_j) > b$ for all $j$.  Thus, the probability that $I$ does not have a local cutpoint is at most $(1-b)^m \leq (1-b)^{C_1 (\log n)^{\theta}}$.  Since there were at most $2C' (\log n)^{1/2 - \theta}$ intervals of length $n^2 (\log n)^{\theta}$, the probability that all of them have a local cutpoint is at least $$1 - 2C' (\log n)^{1/2 - \theta} (1-b)^{C_1 (\log n)^{\theta}},$$ which completes the proof.
\end{proof}

\section{Properties of loop-erased segments}

In this section, we will consider a segment of a loop-erased random walk on the torus $\Z^4_n$.  Our main goals are to show that the length and the capacity are highly concentrated around their means.  As noted in the introduction, we will study the uniform spanning tree on $\Z^4_n$ by breaking up the random walks used in Wilson's algorithm into segments of length $r = \lfloor n^2 (\log n)^{9/22} \rfloor$.  Therefore, in this section we will be primarily concerned with random walks whose length is of order $n^2 (\log n)^{9/22}$.

\subsection{Length of a loop-erased segment}

In this subsection, we show that the length of a segment of a loop-erased random walk on $\Z^4_n$ is highly concentrated around its mean.  The argument relies heavily on results of Lawler \cite{lawler} for the length of a loop-erased random walk on $\Z^4$.  

\begin{Prop}
Suppose $X = (X_t)_{t=0}^{\infty}$ is a random walk on $\Z_n^4$.  Then there exist constants $C$ and $C'$ and a sequence $(c_k)_{k=1}^{\infty}$ with $0 < \inf c_k \leq \sup c_k < \infty$ such that for all $L$,
$$P \bigg( \bigg| \big| LE(X[0,L]) \big| - \frac{L c_L}{(\log L)^{1/3}} \bigg| \geq \frac{C'L}{(\log L)^{3/4}} \bigg) \leq C \bigg(\frac{\log \log L}{(\log L)^{1/4}} + \frac{L^2}{n^4 \log n} \bigg).$$
\label{lengthprop}
\end{Prop}

By applying this proposition when $L$ is of the order $n^2 (\log n)^{9/22}$, we immediately get the following corollary.

\begin{Cor}
Let $X = (X_t)_{t=0}^{\infty}$ be a random walk on $\Z_n^4$.  Assume there are constants $C_1$ and $C_2$ such that $C_1 n^2 (\log n)^{9/22} \leq L \leq C_2 n^2 (\log n)^{9/22}$.  Then there exist constants $C$ and $C'$ and a sequence $(b_n)_{n=1}^{\infty}$ with $0 < \inf b_n \leq \sup b_n < \infty$, all possibly depending on the values of $L$, such that for all $n$, $$P \bigg( \bigg| \big|LE(X[0,L]) \big| - b_n n^2 (\log n)^{5/66} \bigg| \geq \frac{C' n^2}{(\log n)^{15/44}} \bigg) \leq \frac{C}{(\log n)^{4/22}}.$$
\label{lengthcor}
\end{Cor}

\begin{proof}[Proof of Proposition \ref{lengthprop}]
We may assume that $X$ is coupled with a lazy random walk $Z = (Z_t)_{t=0}^{\infty}$ on $\Z^4$ such that $X_t \equiv Z_t$ (mod $n$) for all $t$.  Unless there exist two times $s$ and $t$ such that $X_s = X_t$ and $Z_s \neq Z_t$, we have $|LE(X[0,L])| = |LE(Z[0,L])|$ because the two random walks make loops at the same times.  Therefore, by Proposition \ref{torusintprop},
\begin{equation}
P(|LE(X[0,L])| \neq |LE(Z[0,L])|) \leq \frac{C L^2}{n^4 \log n},
\label{Z4couple}
\end{equation}
so we may focus on controlling $|LE(Z[0,L])|$.

The lazy random walk $Z$ can be coupled with a simple random walk $Y = (Y_t)_{t=0}^{\infty}$ in such a way that if $T_0 = 0$ and $T_k = \min\{t > T_{k-1}: Z_t \neq Z_{T_{k-1}}\}$ for $k \geq 1$, then $Y_t = Z_{T_t}$ for all $t$.  Let $M$ be the cardinality of $\{t: 1 \leq t \leq L \mbox{ and }X_t \neq X_{t-1}\}$, which is independent of $Y[0,L]$ and has the binomial distribution with parameters $L$ and $1/2$.  Then
\begin{equation}
|LE(Z[0,L])| = |LE(Y[0,M])|.
\label{YZcouple}
\end{equation}
It will at times be convenient to condition on $M = m$, and we denote these conditional probabilities and expectations by $P_m$ and $E_m$.

Let $\sigma_0 = \max\{t: Y_t = 0\}$, and for $i \geq 1$, let $\sigma_i = \max\{t > \sigma_{i-1}: Y_t = Y_{\sigma_{i-1} + 1}\}$.  The times $\sigma_0, \sigma_1, \dots$ are said to be retained after loop-erasure.  Let $I_k$ be the indicator of the event that $k = \sigma_i$ for some $i$.  
Let $a_k = E[I_k]$ be the probability that the time $k$ is retained after loop-erasure.  Lawler \cite{law95} showed that there are constants $C_1$ and $C_2$ such that
\begin{equation}
\frac{C_1}{(\log k)^{1/3}} \leq a_k \leq \frac{C_2}{(\log k)^{1/3}}
\label{laweq6}
\end{equation}
for all $k \geq 2$.

Now $\sum_{k=1}^m I_k$ is the number of points, out of the first $m$, that are retained after loop-erasure of the infinite walk.  Call a point $t$ loop-free if, whenever $0 \leq s \leq t$ and $u > t$, we have  $Y_s \neq Y_u$.  By Lemma 7.7.4 of \cite{lawler}, there is a constant $C$ such that for all $m$, the probability that there is a loop-free point between $m - m (\log m)^{-6}$ and $m$ is at least $1 - C (\log \log m)/(\log m)$.  On this event, no points of the walk before time $m - m (\log m)^{-6}$ can be erased after time $m$.  Therefore,
\begin{equation}
P_m \bigg( \bigg| \sum_{k=1}^m I_k - |LE(Y[0,m])| \bigg| > \frac{m}{(\log m)^6} \bigg) \leq \frac{C \log \log m}{\log m}.
\label{laweq3}
\end{equation}

We now consider a finite loop-erasing procedure also considered in \cite{lawler}.  
Choose $0 = j_0 < j_1 < \dots < j_v = M$ such that $M/2(\log M)^2 \leq j_i - j_{i-1} \leq 2M/(\log M)^2$ for $i = 1, \dots, v$.  If $j_{i-1} < k \leq j_i$, or if $i = 1$ and $k = 0$, let $I'_k$ be the indicator of the event that the time $k$ is retained after loop-erasure of the segment $X[j_{i-1}, j_i]$.  Lawler shows (see p. 207 of \cite{lawler}) that
\begin{equation}
E_m \bigg[ \bigg| \sum_{k=0}^m I_k - \sum_{k=0}^m I'_k \bigg| \bigg]  \leq \frac{C m \log \log m}{\log m}.
\label{laweq2}
\end{equation}
Lawler defines the times $j_i$ slightly differently, but it is easy to see that this does not affect the proof of (\ref{laweq2}).  Therefore, by Markov's Inequality,
\begin{equation}
P_m \bigg( \bigg| \sum_{k=0}^m I_k - \sum_{k=0}^m I'_k \bigg| \geq \frac{m}{(\log m)^{3/4}} \bigg) \leq \frac{C \log \log m}{(\log m)^{1/4}}. 
\label{laweq4}
\end{equation}
Lawler also shows (see p. 208 of \cite{lawler}) that
\begin{equation}
P_m \bigg( \bigg| \sum_{k=0}^m I'_k - E_m \bigg[ \sum_{k=0}^m I'_k \bigg] \bigg| >
\frac{1}{(\log m)^{1/2}} E_m \bigg[ \sum_{k=0}^m I'_k \bigg] \bigg) \leq \frac{C}{(\log m)^{3/8}}.
\label{laweq7}
\end{equation}
Since $E[\sum_{k=1}^m I_k] = \sum_{k=1}^m a_k \leq Cm/(\log m)^{1/3}$ by (\ref{laweq6}), we can combine (\ref{laweq7}) with (\ref{laweq2}) to get that there are constants $C$ and $C'$ such that for all $m$,
\begin{equation}
P_m \bigg( \bigg| \sum_{k=0}^m I'_k - E_m \bigg[ \sum_{k=0}^m I'_k \bigg] \bigg| >
\frac{C' m}{(\log m)^{5/6}} \bigg) \leq \frac{C}{(\log m)^{3/8}}.
\label{laweq1}
\end{equation}
By (\ref{laweq3}), (\ref{laweq4}), and (\ref{laweq1}), there are constants $C$ and $C'$ such that
\begin{equation}
P_m \bigg( \bigg| |LE((Y_t)_{t=0}^m)| - E_m \bigg[ \sum_{k=0}^m I'_k \bigg] \bigg| > \frac{C' m}{(\log m)^{3/4}} \bigg) \leq \frac{C \log \log m}{(\log m)^{1/4}}.
\label{laweq5}
\end{equation}

By (\ref{azuma}),
\begin{equation}
P(|M - L/2| > L^{3/4}) \leq 2 e^{-2 \sqrt{L}}.
\label{LDML}
\end{equation}
Therefore,
\begin{equation}
P \bigg( \bigg| \sum_{k=0}^M a_k - \sum_{k=0}^{\lfloor L/2 \rfloor} a_k \bigg| > 1 + L^{3/4} \bigg) \leq 2 e^{-2 \sqrt{L}}.
\label{laweq8}
\end{equation}
We have
\begin{align}
\bigg| \big| LE&(Y[0,M]) \big| - \sum_{k=0}^{\lfloor L/2 \rfloor} a_k \bigg| \nonumber \\
&\leq \bigg| \big| LE(Y[0,M]) \big| -
E \big[ \sum_{k=0}^M I_k' \big| M \big] \bigg| + \bigg| E \big[ \sum_{k=0}^M I_k' \big| M \big] - E \big[ \sum_{k=0}^M  I_k \big| M \big] \bigg| + \bigg| \sum_{k=0}^M a_k - \sum_{k=0}^{\lfloor L/2 \rfloor} a_k \bigg|. \nonumber
\end{align}
By conditioning on the value of $M$, we can use
(\ref{laweq5}) and (\ref{LDML}) to bound the first term, (\ref{laweq2}) and (\ref{LDML}) to bound the second term, and (\ref{laweq8}) to bound the third term, we get
\begin{equation}
P \bigg( \bigg| \big| LE(Y[0,M]) \big| - \sum_{k=0}^{\lfloor L/2 \rfloor} a_k \bigg| > \frac{C'L}{(\log L)^{3/4}} \bigg) \leq \frac{C \log \log L}{(\log L)^{1/4}}.
\label{laweq9}
\end{equation}
For all $k$, let $$c_k = \frac{(\log k)^{1/3}}{k} \sum_{j=0}^{\lfloor k/2 \rfloor} a_j.$$  It is easily seen from (\ref{laweq6}) that $0 < \inf c_k < \sup c_k < \infty$.  The proposition now follows from (\ref{laweq9}), (\ref{Z4couple}), and (\ref{YZcouple}).
\end{proof}

\subsection{Capacity of a loop-erased segment}

Let $U$ and $V$ be subsets of $\Z_n^4$.  Let $(X_t)_{t=0}^{\infty}$ be a random walk on $\Z_n^4$ such that $X_0$ has the uniform distribution on $\Z_n^4$.  Recall from (\ref{capdef}) that the capacity of a set $U \subset \Z^4_n$ is defined by
$$\mbox{Cap}_M(U) = P(X_t \in U \mbox{ for some }t \leq M).$$  As in \cite{perrev}, define the closeness of the sets $U$ and $V$ by $$\mbox{Close}_M(U, V) = P(X_t \in U \mbox{ and }X_s \in V \mbox{ for some }s, t \leq M).$$  In this subsection, we estimate the capacity of a loop-erased random walk segment, and show that the capacity is highly concentrated around its mean.  The first step is the following bound on the closeness between two random walk segments.

\begin{Prop}
Let $X = (X_t)_{t=0}^{\infty}$ and $Y = (Y_t)_{t=0}^{\infty}$ be independent random walks on $\Z^4_n$, and assume that $Y_0$ has the uniform distribution on $\Z^4_n$.  Assume there is a constant $C'$ such that $K \geq C' n \log n$, $L \geq C' n \log n$, $M \geq C' n \log n$, and $u \geq \tau_n$.  Let $U = X[0,K]$, $V = Y[0,L]$, and $W = X[K+u, K+L+u]$.  Then there is a constant $C$ such that
$$E[\textup{Close}_M(U, V)] \leq \frac{C K L M^2}{n^8 (\log n)^2} \hspace{.2in}\mbox{and}\hspace{.2in}
E[\textup{Close}_M(U, W)] \leq \frac{C K L M^2}{n^8 (\log n)^2} .$$
\label{closeprop}
\end{Prop}

\begin{proof}
The proof is similar to the proof of Lemma 5.2 in \cite{perrev}.  We may assume that $X_0$ has the uniform distribution on $\Z^4_n$.
Let $Z = (Z_t)_{t=0}^{\infty}$ be a random walk on $\Z^4_n$, independent of $X$ and $Y$, such that $Z_0$ has the uniform distribution on $\Z^4_n$.  For all sets $S$, let $T_S = \min\{t: Z_t \in S\}$.  We have
\begin{align}
E[\textup{Close}_M(U, V)] &\leq P(T_U \leq M) P(Z_t \in V \mbox{ for some }T_U \leq t \leq M|T_U \leq M) \nonumber \\
&\hspace{.5in}+ P(T_V \leq M) P(Z_t \in U \mbox{ for some }T_V \leq t \leq M|T_V \leq M),
\label{2termclose}
\end{align}
and the same inequality holds with $W$ in place of $V$.  By Proposition \ref{mainintprop}, there is a constant $C$ such that
\begin{align}
P(T_U \leq M) &= E[\textup{Cap}_M(U)] \leq \frac{C K M}{n^4 \log n}, \nonumber \\
P(T_V \leq M) &= E[\mbox{Cap}_M(V)] \leq \frac{C L M}{n^4 \log n}, \nonumber \\
P(T_W \leq M) &= E[\mbox{Cap}_M(W)] \leq \frac{C L M}{n^4 \log n}. \nonumber
\end{align}
The strong Markov property applied at $T_U$ implies that
$$P(Z_t \in V \mbox{ for some }T_U \leq t \leq M|T_U \leq M) \leq E[\mbox{Cap}_M(V)] \leq \frac{CLM}{n^4 \log n}.$$  Likewise, the strong Markov property applied at $T_V$ gives
$$P(Z_t \in U \mbox{ for some }T_V \leq t \leq M|T_V \leq M) \leq E[\mbox{Cap}_M(U)] \leq \frac{CKM}{n^4 \log n}.$$  It follows that the right-hand side of (\ref{2termclose}) is at most $2C^2 KLM^2/[n^8 (\log n)^2]$.  Since $u \geq \tau_n$, the bound for $E[\mbox{Close}_M(U,W)]$ now follows from Lemma \ref{mixlem} with $f(U, V) = \mbox{Close}_M(U,V)$.
\end{proof}

\begin{Lemma}
Let $Y = (Y_t)_{t=0}^{\infty}$ be a random walk on $\Z^4_n$.  Let $U = \{x: Y_t = x \mbox{ for some }t \leq L\}$.  Assume that $L \geq C' n \log n$ and $M \geq C' n \log n$ for some constant $C'$.
Then there exists a constant $C$ such that $$E[\textup{Cap}_M(U)^2] \leq \frac{C L^2 M^2}{n^8 (\log n)^2}.$$
\label{varcaplemma}
\end{Lemma}

\begin{proof}
Let $X = (X_t)_{t=0}^{\infty}$ and $X' = (X'_t)_{t=0}^{\infty}$ be two independent random walks on $\Z_n^4$ such that $X_0$ and $X'_0$ have the uniform distribution.  Then
\begin{align}
E[\mbox{Cap}_M(U)^2] &= E[P(X[0,M] \cap U \neq \emptyset|U)^2] \nonumber \\
&= E[P(X[0,M] \cap U \neq \emptyset|U) P(X'[0,M] \cap U \neq \emptyset|U)] \nonumber \\
&= E[P(X[0,M] \cap U \neq \emptyset, X'[0,M] \cap U \neq \emptyset|U)] \nonumber \\
&= E[P(X[0,M] \cap U \neq \emptyset, X'[0,M] \cap U \neq \emptyset|X, X')] \nonumber \\
&= E[\mbox{Close}_L(X[0,M], X'[0,M])] \nonumber \\
&\leq \frac{C L^2 M^2}{n^8 (\log n)^2}
\end{align}
by Proposition \ref{closeprop}.
\end{proof}

We now work towards a lower-bound on the probability that a segment of a random walk on $\Z^4_n$ intersects a segment of a loop-erased random walk on $\Z^4_n$.  Lyons, Peres, and Schramm \cite{lps03} showed that for two transient Markov chains with the same transition probabilities, the probability that one intersects the loop-erasure of the other is at least $2^{-8}$ times the probability that the paths of the two chains intersect.  That is, erasing loops for one of the two chains reduces the intersection probability only by a constant.  Because this result requires the Markov chains to be transient, it applies to random walks on $\Z_n^4$ only if the walks are killed at a geometric time.  In Proposition \ref{lpsprop} below, we adapt this result to obtain a lower bound on the intersection probability for walks of fixed length.  The proof uses the following lemma pertaining to the geometric distribution.

\begin{Lemma}
Suppose $P(X = k) = p (1-p)^{k-1}$ for $k = 1, 2, \dots$.  Then for all positive integers $m$, we have $E[X {\bf 1}_{\{X > m\}}] = (1-p)^m (m + 1/p)$.
\label{geomlem}
\end{Lemma}

\begin{proof}
We have
\begin{align}
E[X {\bf 1}_{\{X > m\}}] &= \sum_{k=1}^{\infty} P(X {\bf 1}_{\{X > m\}} \geq k)
= m P(X > m) + \sum_{k=m+1}^{\infty} P(X \geq k) \nonumber \\
&= m (1-p)^m + \sum_{k=m+1}^{\infty} (1-p)^{k-1}  = (1-p)^m \bigg( m + \frac{1}{p} \bigg), \nonumber
\end{align}
as claimed.
\end{proof}

\begin{Prop}
Let $X = (X_t)_{t=0}^{\infty}$ and $Y = (Y_t)_{t=0}^{\infty}$ be independent random walks on $\Z^4_n$.  Assume that $Y_0$ has the uniform distribution on $\Z_n^4$.  If $L/n^2 \rightarrow \infty$ and $L/(n^2 (\log n)^{1/2}) \rightarrow 0$ as $n \rightarrow \infty$, then there is a positive constant $C$ such that $$P(LE(X[0,L]) \cap Y[0,L] \neq \emptyset) \geq \frac{C L^2}{n^4 \log n}.$$
\label{lpsprop}
\end{Prop}

\begin{proof}
Let $A \geq 1$ be a number, to be chosen more precisely later, such that $L/A$ is an integer greater than one.  Let $S$ and $T$ be independent random variables, independent of $X$ and $Y$, having a geometric distribution with parameter $p = A/L$.  Therefore, $P(S > L/A) = P(T > L/A) = (1 - A/L)^{L/A} \geq 1/4$.
By Proposition \ref{lowerintprop}, there is a constant $C_1$ such that
\begin{align}
P(X[0,S] \cap Y[0,T] \neq \emptyset) &\geq P(S > L/A, T > L/A, X[0, L/A] \cap Y[0, L/A] \neq \emptyset) \nonumber \\
&= P(S > L/A) P(T > L/A) P(X[0, L/A] \cap Y[0, L/A] \neq \emptyset) \nonumber \\
&\geq \frac{C_1 L^2}{A^2 n^4 \log n}. \nonumber
\end{align}
By Lemma 1.2 of \cite{lps03},
\begin{equation}
P(LE(X[0,S]) \cap Y[0,T] \neq \emptyset) \geq 2^{-8} P(X[0,S] \cap Y[0,T] \neq \emptyset) \geq \frac{C_1 L^2}{2^{8} A^2 n^4 \log n}.
\label{lpseq}
\end{equation}

The goal now is to choose $A$ so that if $LE(X[0,S])$ and $Y[0,T]$ intersect, then $LE(X[0,L])$ and $Y[0,L]$ also intersect outside an event of low probability.  We define the following events:
\begin{itemize}
\item Let $B_1$ be the event that $LE(X[0,S]) \cap Y[0,T] \neq \emptyset$ and $\max\{S, T\} > L$.

\item Let $B_2$ be the event that $X[0, S - \tau_n] \cap X[S, L] \neq \emptyset$ and $X[0, S - \tau_n] \cap Y[0,L] \neq \emptyset$.

\item Let $B_3$ be the event that $X[S - \tau_n, S] \cap Y[0, L] \neq \emptyset$.
\end{itemize}
Note that if $LE(X[0,S]) \cap Y[0,T] \neq \emptyset$ but $LE(X[0,L]) \cap Y[0,L] = \emptyset$, then $B_1 \cup B_2 \cup B_3$ must occur.  The event $B_1$ accounts for the possibility that the segments $LE(X[0,L])$ and $Y[0,L]$ could fail to intersect because the walks are run for a shorter time.  If $L \geq \max\{S, T\}$, then the segments coming from longer walks could still fail to intersect if the portion of the path $LE(X[0,S])$ that $Y[0, T]$ hits gets erased between times $S$ and $L$.  If this happens then either $B_2$ or $B_3$ must occur, with $B_3$ accounting for the possibility that $Y$ hits near the end of the path $LE(X[0,S])$ and $B_2$ accounting for the possibility that $Y$ hits earlier in the path but this part of the path has been erased by an intersection with a later part of the $X$ path.  Thus, to obtain the proposition from (\ref{lpseq}), we must bound $P(B_1)$, $P(B_2)$, and $P(B_3)$.

To bound $P(B_1)$, note that by Proposition \ref{mainintprop}, there is a constant $C_2$ such that on the event $S \geq n \log n$ and $T \geq n \log n$,
$$P(X[0,S] \cap Y[0,T] \neq \emptyset|S, T) \leq \frac{C_2 ST}{n^4 \log n}.$$
Therefore, by conditioning on $S$ and $T$, we get
\begin{align}
P(B_1) &\leq P(X[0,S] \cap Y[0,T] \neq \emptyset \mbox{ and } \max\{S, T\} > L) \nonumber \\
&\leq 2 P(X[0,S] \cap Y[0,T] \neq \emptyset \mbox{ and } S > L) \nonumber \\
&\leq \frac{2 C_2}{n^4 \log n} E[S T {\bf 1}_{\{ S > L\}}] + 2 P(T < n \log n). \nonumber
\end{align}
We know that $S$ and $T$ are independent and $E[T] = 1/p = L/A$.  By Lemma \ref{geomlem},
$$E[S {\bf 1}_{\{S > L\}}] = \bigg( 1 - \frac{A}{L} \bigg)^L \bigg( L + \frac{L}{A} \bigg)
\leq e^{-A} \bigg( L + \frac{L}{A} \bigg) \leq 2 e^{-A} L.$$
Also, $P(T < n \log n) \leq (A n \log n)/L$.  It follows that
\begin{equation}
P(B_1) \leq \frac{4 C_2 e^{-A} L^2}{A n^4 \log n} + \frac{2A n \log n}{L}.
\label{PB1}
\end{equation}

To bound $P(B_2)$, note that conditional on $X[0, S - \tau_n]$, the probability that $X[0, S - \tau_n] \cap Y[0,L] \neq \emptyset$ is $\mbox{Cap}_L(X[0, S - \tau_n])$, and by Lemma \ref{mixlem} the probability that $X[0, S - \tau_n] \cap X[S, L] \neq \emptyset$ is at most $(3/2) \mbox{Cap}_L(X[0, S - \tau_n])$.  Because $X[S, L]$ and $Y[0,L]$ are conditionally independent given $X[0, S - \tau_n]$, we have
\begin{equation}
\label{PB2}
P(B_2) \leq \frac{3}{2}E[\mbox{Cap}_L(X[0, S - \tau_n])^2] \leq \frac{3}{2} E[\mbox{Cap}_L(X[0,S])^2] \leq \frac{CL^4}{A^2 n^8 (\log n)^2},
\end{equation}
where the last inequality follows from Lemma \ref{varcaplemma} by conditioning on $S$.  Also by Proposition \ref{mainintprop},
\begin{equation}
\label{PB3}
P(B_3) \leq \frac{C \tau_n L}{n^4 \log n}.
\end{equation}

By choosing $A$ sufficiently large, we can make the right-hand side of (\ref{PB1}) smaller than one-half the right-hand side of (\ref{lpseq}).  The right-hand of (\ref{PB2}) is much smaller than the right-hand side of (\ref{lpseq}) because $L/(n^2 (\log n)^{1/2}) \rightarrow 0$, and the right-hand side of (\ref{PB3}) is much smaller than the right hand side of (\ref{lpseq}) because $L/\tau_n \rightarrow \infty$ in view of (\ref{mixeq}).  The proposition follows.
\end{proof}

We now use Proposition \ref{closeprop} to bound the probability that a segment of a random walk intersects a shorter loop-erased segment.  The result does not follow immediately from Proposition \ref{lpsprop} because in Proposition \ref{lpsprop}, we require the two segments to be the same length.  This assumption was required for the proof because, to apply the result of Lyons, Peres, and Schramm \cite{lps03}, we needed to assume that the two killed random walks had the the same distribution, and therefore that they were killed at the same rate.  However, we can extend the result to cover the case in which the loop-erased segment is shorter by splitting the longer random walk segment into pieces.

\begin{Prop}
Let $(X_t)_{t=0}^{\infty}$ and $(Y_t)_{t=0}^{\infty}$ be independent random walks on $\Z^4_n$ such that $Y_0$ has the uniform distribution on $\Z^4_n$.  Assume that $K \leq L$ with $K/n^2 \rightarrow \infty$ and $L/(n^2 (\log n)^{1/2}) \rightarrow 0$ as $n \rightarrow \infty$.  Then there exist constants $C_1$ and $C_2$ such that
\begin{equation}
\frac{C_1 K L}{n^4 \log n} \leq P(LE(X[0,K]) \cap Y[0,L] \neq \emptyset) \leq \frac{C_2 K L}{n^4 \log n}.
\label{lpspropeq}
\end{equation}
\label{lpssmallprop}
\end{Prop}

\begin{proof}
We may assume that $X_0$ also has the uniform distribution on $\Z^4_n$.
Because $$P(LE(X[0,K]) \cap Y[0,L] \neq \emptyset) \leq
P(X[0,K] \cap Y[0,L] \neq \emptyset),$$ the upper bound in the proposition follows immediately from Proposition \ref{mainintprop}.

To prove the lower bound, let $M = \lfloor L/2K \rfloor$.  If $M = 0$, then $L < 2K$ and because
\begin{equation}
P(LE(X[0,K]) \cap Y[0,L] \neq \emptyset) \geq P(LE(X[0,K]) \cap Y[0,K] \neq \emptyset),
\label{LEKL}
\end{equation}
the lower bound follows from Proposition \ref{lpsprop}.  Assume now that $M \geq 1$.  For $j = 0, 1, \dots, M-1$, we have $$P(LE(X[0,K]) \cap Y[2jK, (2j+1)K] \neq \emptyset) \geq \frac{C_3 K^2}{n^4 \log n}$$ for some positive constant $C_3$ by Proposition \ref{lpsprop}.  The probability that $LE(X[0,K])$ intersects both $Y[2jK, (2j+1)K]$ and $Y[2 \ell K, (2 \ell + 1)K]$ is at most the probability that $X[0,K]$ intersects both of these segments, which is $E[\textup{Close}_K(Y[2jK, (2j+1)K], Y[2 \ell K, (2 \ell + 1)K])]$.  If $j \neq \ell$, the intervals $[2jK, (2j+1)K]$ and $[2 \ell K, (2 \ell + 1)K]$ are separated by a gap of at least $K \geq \tau_n$, so by Proposition \ref{closeprop}, $$E[\textup{Close}_K(Y[2jK, (2j+1)K], Y[2 \ell K, (2 \ell + 1)K])] \leq \frac{C_4 K^4}{n^8 (\log n)^2}$$ for some constant $C_4$.  There are $M$ of these segments and $M(M-1)/2$ pairs of segments, so by inclusion-exclusion,
$$P(LE(X[0,K]) \cap Y[0,L] \neq \emptyset) \geq \frac{C_3 M K^2}{n^4 \log n} - \frac{C_4 M (M-1) K^4}{2 n^8 (\log n)^2} \geq \frac{C_1 KL}{n^4 \log n}$$ for some constant $C_1$, as claimed.
\end{proof}

\begin{Prop}
Let $X = (X_t)_{t=0}^{\infty}$ be a random walk on $\Z^4_n$.  Assume that there are constants $C_1$ and $C_2$ such that $C_1 n^2 (\log n)^{9/22} \leq L \leq C_2 n^2 (\log n)^{9/22}$ and $C_1 n^2 (\log n)^{9/22} \leq M \leq C_2 n^2 (\log n)^{9/22}$.  Then, for all $\theta > 0$, there is a constant $C$ and a sequence of constants $(a_n)_{n=1}^{\infty}$ with $0 < \inf a_n \leq \sup a_n < \infty$, possibly depending on the values of $L$ and $M$, such that
$$P \bigg( \bigg| \textup{Cap}_M(LE(X[0,L])) - \frac{a_n}{(\log n)^{2/11}} \bigg| > \frac{1}{(\log n)^{5/22}} \bigg) \leq \frac{C}{(\log n)^{3/22 - \theta}}.$$
\label{maincapprop}
\end{Prop}

\begin{proof}
Let $k_j = j \lceil n^2 (\log n)^{4/22} \rceil$, and define $v$ so that $0 = k_0 < k_1 < \dots < k_v \leq L < k_{v+1}$.  Note that $v \leq L/(n^2 (\log n)^{4/22}) \leq C_2 (\log n)^{5/22}$.  Let $U = LE(X[0,L])$, and for $i = 1, \dots, v$, let $V_i = LE(X[k_{i-1}, k_i])$ and $W_i = X[k_{i-1}, k_i]$.  Let $D$ be the set of points in $U$ or in $\cup_{i=1}^v V_i$ but not both.  We first claim that
\begin{equation}
\sum_{i=1}^v \mbox{Cap}_M(V_i) - \mbox{Cap}_M(D) - \sum_{i=1}^v \sum_{j=1}^{i-1} \mbox{Close}_M(V_i, V_j) \leq \mbox{Cap}_M(U) \leq \sum_{i=1}^v \mbox{Cap}_M(V_i) + \mbox{Cap}_M(D).
\label{capclose}
\end{equation}
To see this, let $Y = (Y_t)_{t=0}^{\infty}$ be a random walk on $\Z^4_n$ started from the uniform distribution.  If the path $Y[0,M]$ intersects $U$, then it must either intersect $D$ or one of the $V_i$, which gives the upper bound for $\mbox{Cap}_M(U)$.  For the lower bound, first note that, conditional on $X[0,L]$, the probability that $Y[0,M]$ intersects $U$ but not $\cup_{i=1}^v V_i$ is at most the probability that it intersects $D$.  By inclusion-exclusion, the probability that it intersects $\cup_{i=1}^v V_i$ is at most the sum of the probabilities that it intersects the individual $V_i$ minus the sum over pairs $i$ and $j$ of the probabilities that it intersects both $V_i$ and $V_j$.

Let $B_1$ be the event that there exist $s$ and $t$ with $s + 2 \tau_n \leq t \leq L$ such that $X_s = X_t$.  That is, $B_1$ is the event that the path $X[0,L]$ has a long loop of length at least $2 \tau_n$.  Fix $\theta > 0$, and let $w = \lfloor n^2 (\log n)^{\theta} \rfloor$.  Let $B_2$ be the event that the path $X[0,L]$ has an interval $[t, t + w]$ with $0 \leq t \leq L - n^2 (\log n)^{\theta}$ that does not contain a local cutpoint.  By Proposition \ref{medloops}, $P(B_1) \leq C (\log n)^{-4/22}$, and by Proposition \ref{lcprop}, $P(B_2) \leq C (\log n)^{-1}$.

Let $B = (B_1 \cup B_2)^c$.  On the event $B$, for each $i = 1, 2, \dots, v$, the path $X[0,L]$ has local cutpoints $u_i$ and $v_i$ such that $u_i \in [k_{i-1}, k_{i-1} + w]$ and $v_i \in [k_i - w, k_i]$.  On the event $B$, if $t_1 < u_i < t_2 < v_i < t_3$, then $X_{t_1}$, $X_{t_2}$, and $X_{t_3}$ must all be distinct, as the definition of a local cutpoint prohibits short loops while long loops do not occur on $B_1^c$.  It follows that if $u_i \leq t \leq v_i$, then $X_t$ is a point in the path $LE(X[k_{i-1}, k_i - 1])$ if and only if $X_t$ is a point in the path $LE(X[0,L])$.  That is, we have $X_t \in U$ if and only if $X_t \in V_i$.  It follows that $D \subset D'$, where $D' = X[k_{i-1}, k_{i-1}+w] \cup X[k_i - w, k_i] \cup X[k_v, L]$.
By Proposition \ref{mainintprop}, the expected capacity of a random walk segment of length $w + 1$ is at most $C M w/(n^4 \log n)$.  Because the capacity of $D'$ is at most the sum of the capacities of $2v$ segments of length $w + 1$ plus the capacity of one segment of length at most $n^2 (\log n)^{4/22}$, it follows that
\begin{align}
E[\mbox{Cap}_M(D) {\bf 1}_{B}] &\leq E[\mbox{Cap}_M(D')] \leq 2v \bigg( \frac{C M w}{n^4 \log n} \bigg) + \frac{CM n^2 (\log n)^{4/22}}{n^4 \log n} \nonumber \\
&\leq C (\log n)^{5/22 + 9/22 + \theta - 1} + C (\log n)^{9/22 + 4/22 - 1} \leq C (\log n)^{-8/22 + \theta}.
\label{caprbound}
\end{align}

We next bound $\mbox{Close}_M(V_i, V_j)$.  We have $\mbox{Close}_M(V_i, V_j) \leq \mbox{Close}_M(W_i, W_j)$, as erasing loops can only reduce the probability of intersection.  Therefore, by Proposition \ref{closeprop}, if $j < i - 1$ there is a constant $C$ such that
$$E[\mbox{Close}_M(V_i, V_j)] \leq \frac{C M^2 (k_{i} - k_{i-1})(k_{j} - k_{j-1})}{n^8 (\log n)^2}
\leq \frac{4 C n^8 (\log n)^{18/22 + 8/22}}{n^8 (\log n)^2} \leq 4C (\log n)^{-18/22}.$$
To bound the closeness between $V_i$ and $V_{i-1}$, note that for the random walk $Y[0,M]$ to intersect both $V_i$ and $V_{i-1}$, it must either intersect $X[k_{i-1}-w, k_{i-1}]$, which has probability at most $C M w/(n^4 \log n) \leq C (\log n)^{-13/22 + \theta}$ by Proposition \ref{mainintprop}, or else it must intersect both $X[k_{i-2}, k_{i-1} - w]$ and $X[k_{i-1}, k_i - 1]$, which has probability at most $C (\log n)^{-18/22}$ by Proposition \ref{closeprop}.  Therefore, $E[\mbox{Close}_M(V_{i-1}, V_i)] \leq C (\log n)^{-13/22 + \theta}$.  It follows that
\begin{equation}
E \bigg[ \sum_{i=1}^v \sum_{j=1}^{i-1} \mbox{Close}_M(V_i, V_j)  \bigg] \leq Cv (\log n)^{-13/22 + \theta} + C \binom{v}{2} (\log n)^{-18/22} \leq C (\log n)^{-8/22 + \theta}.
\label{closerbound}
\end{equation}
By combining (\ref{capclose}), (\ref{caprbound}), and (\ref{closerbound}), we get
$$E \bigg[ \bigg| \mbox{Cap}_M(U) - \sum_{i=1}^v \mbox{Cap}_M(V_i) \bigg| {\bf 1}_B \bigg] \leq C (\log n)^{-8/22 + \theta}.$$
Now Markov's Inequality gives
\begin{equation}
P \bigg( \bigg| \mbox{Cap}_M(U) - \sum_{i=1}^v \mbox{Cap}_M(V_i) \bigg| > \frac{1}{2(\log n)^{5/22}} \bigg) \leq C (\log n)^{-3/22 + \theta} + P(B^c) \leq C (\log n)^{-3/22 + \theta}.
\label{Markbound}
\end{equation}

By Lemma \ref{varcaplemma}, we have
$$\mbox{Var}(\mbox{Cap}_M(V_i)) \leq E[(\mbox{Cap}_M(V_i))^2] \leq E[(\mbox{Cap}_M(W_i))^2] \leq C (\log n)^{-18/22}.$$  Because the random variables $\mbox{Cap}_M(V_1), \dots, \mbox{Cap}_M(V_v)$ are independent, it follows that
$$\mbox{Var} \bigg( \sum_{i=1}^v \mbox{Cap}_M(V_i) \bigg) \leq C (\log n)^{-13/22}.$$  Therefore, by Chebyshev's Inequality,
\begin{equation}
P \bigg( \bigg| \sum_{i=1}^v \mbox{Cap}_M(V_i) - E \bigg[ \sum_{i=1}^v \mbox{Cap}_M(V_i) \bigg] \bigg| > \frac{1}{2(\log n)^{5/22}} \bigg) \leq C (\log n)^{-3/22}.
\label{Chebbound}
\end{equation}

By Proposition \ref{lpssmallprop}, there exist constants $C_1$ and $C_2$ such that for all $i = 1, \dots, v$, we have $C_1 (\log n)^{-9/22} \leq E[\mbox{Cap}_M(V_i)] \leq C_2 (\log n)^{-9/22}$, and therefore there exist constants $C_3$ and $C_4$ such that $$C_3 (\log n)^{-2/11} \leq E \bigg[ \sum_{i=1}^v \mbox{Cap}_M(V_i) \bigg] \leq C_4 (\log n)^{-2/11}.$$  Therefore, by setting $a_n = (\log n)^{2/11} E[ \sum_{i=1}^v \mbox{Cap}_M(V_i)]$, we get the result from (\ref{Markbound}) and (\ref{Chebbound}).
\end{proof}

\subsection{Decomposing and labeling a loop-erased path}

Let $X = (X_t)_{t=0}^{\infty}$ be a random walk on $\Z^4_n$.  Let $T$ be a fixed or random time.  In this subsection, we will consider the loop-erased path $LE(X[0,T])$.  We will study the structure of the loop-erased walk by breaking the original path into segments of length $r = \lfloor n^2 (\log n)^{9/22} \rfloor,$ and keeping track of which of these segments are retained after loop-erasure.  This will lay the groundwork for the coupling in the next section, which will require a similar decomposition of the loop-erased paths used in Wilson's algorithm.

To define the time indices that are retained after loop-erasure, for any positive integers $u \leq v$ let $\sigma_0^{u,v} = \max\{t: u \leq t \leq v \mbox{ and }X_t = X_u\}$.  For $i \geq 1$, whenever $\sigma_{i-1}^{u,v} < v$, let $$\sigma_i^{u,v} = \max\{t: \sigma_{i-1}^{u,v} < t \leq v \mbox{ and }X_t = X_{\sigma_{i-1}^{u,v} + 1}\}.$$
Let $W(u,v)$ denote the set of all times $\sigma_i^{u,v}$, so the path $LE(X[u,v])$ is the same as the path $(X_t)_{t \in W(u,v)}$.  Note that if $v_1 \leq v_2$, then $(W(u, v_2) \cap [u, v_1]) \subset W(u, v_1)$ because new time indices can get erased but not added when the random walk is extended from time $v_1$ to time $v_2$.  We will consider segments of the random walk of length $r$.  Define the $j$th such segment to be $A_j = \{(j-1)r, \dots, jr - 1\}$, and let $\ell = \lfloor T/r \rfloor + 1$.  Let $\theta > 0$, and then set $w = \lfloor n^2 (\log n)^{\theta} \rfloor$.  Let $A_j' = \{(j-1)r + w, \dots, jr - w - 1\}$, which is $A_j$ with the first and last $w$ points removed.  We assume that $n$ is large enough that $w \geq 2 \tau_n$.  Denote $(X_t)_{t \in A_j}$ by $X[A_j]$.

We now define the events $B_1, \dots, B_5$.  When none of them occurs, we say the walk $X[0,T]$ is ``good."  This condition is close to the condition introduced by Peres and Revelle \cite{perrev} for a path to be ``locally decomposable."  When $X[0,T]$ is good, it will be possible to approximate $LE(X[0,T])$ using the loop-erasures of the individual segments of length $r$. 

\begin{itemize}
\item Let $B_1$ be the event that $X[0,T]$ has an interval $[t, t + w]$ with $0 \leq t \leq T - w$ that does not contain a local cutpoint.

\item Let $B_2$ be the event that for some $j = 1, \dots, \ell - 1$, there exist two times $s, t \in A_j \cup A_{j+1}$ such that $t \geq s + 2 \tau_n$ and $X_s = X_t$.

\item Let $B_3$ be the event that there exist distinct integers $i, j, k \leq \ell$ such that $|i-j| \geq 2$, $|i-k| \geq 2$, and $X[A_i]$ intersects both $X[A_j]$ and $X[A_k]$.

\item Let $B_4$ be the event that there exist $i \leq \ell$ and $j \leq \ell$ with $|i - j| \geq 2$ such that $X_s = X_t$ for some $s \in A_i$ and $t \in A_j \setminus A_j'$.

\item Let $B_5$ be the event that for some $j \leq \ell - 2$, we have $X[A_j] \cap X[A_{\ell}] \neq \emptyset$.
\end{itemize}

On $B_1^c$, the walk contains local cutpoints in every segment of length $w$.
On $B_2^c$, there are no long loops within segments of length $2r$.  The event $B_3^c$ prohibits one segment of length $r$ from intersecting two other segments, while $B_4^c$ prevents intersections involving the first or last $w$ steps of one of the segments.  On $B_5^c$, the last segment does not intersect the others.  Later, we will prove Proposition \ref{goodprop}, which will establish that if $T$ has a geometric distribution with a mean of order $n^2 (\log n)^{1/2}$, then the probability that $X[0,T]$ is good is at least $1 - C(\log n)^{-1/11}$.

To keep track of intersections between these segments of length $r$, let $I_{j, j+1} = 0$ for $j = 1, 2, \dots, \ell - 1$.  For integers $i$ and $j$ with $1 \leq i \leq \ell$, $1 \leq j \leq \ell$, and $j \geq i + 2$, define $I_{i,j}$ to be the indicator of the event that $LE(X[A_i]) \cap X[A_j] \neq \emptyset$.  Following \cite{perrev}, let $S_0 = \{0\}$ and, for $j = 1, \dots, \ell-1$, define
$$S_j = \{k \in S_{j-1}: I_{i,j} = 0 \mbox{ for all }i \in \{1, \dots, k\} \cap S_{j-1}\} \cup \{j\}.$$  Roughly speaking, $S_j$ consists of the segments that have survived loop-erasure up to time $jr$.  Let $S = S_{\ell-1}$.  Let $J_j = \{k \leq j: I_{i,k} = 1 \mbox{ or }I_{k,i} = 1 \mbox{ for some }i \leq j\}$ consist of the indices of the segments that are involved in an intersection by time $jr$, and let $J = J_{\ell - 1} \cup \{\ell\}$.  

The next lemma will enable us to approximate the loop-erasure of the whole walk $X[0,T]$ by the loop-erasure of the individual segments when the walk is good.  The lemma states that for $k \in S \cap J^c$, the times retained in $LE(X[A_k])$ are exactly the times in $A_k$ retained in $LE(X[0, jr-1])$, except possibly for the first and last $w$ steps, while no times are retained for $k \in S^c \cap J^c$.

\begin{Lemma}
Suppose $X[0,T]$ is good.  Then $$W(0, T) \cap A_k' = W((k-1)r, kr-1) \cap A_k'$$ for all $k \in S \cap J^c$, and $W(0, T) \cap A_k = \emptyset$ for all $k \in S^c \cap J^c$.
\label{goodlem}  
\end{Lemma}

\begin{proof}
We claim that for all $j = 1, \dots, \ell - 1$, we have
\begin{equation}
W(0, jr-1) \cap A_k' = W((k-1)r, kr-1) \cap A_k'
\label{Weq}
\end{equation}
for $k \in S_j \cap J_j^c$, and $W(0, jr - 1) \cap A_k = \emptyset$ for all $k \in S_j^c \cap J_j^c$. 
To see how the lemma follows from the claim, note that on $B_5^c$, we have $W(0, (\ell-1)r - 1) \cap A_k' = W(0, T) \cap A_k'$ for all $k \leq \ell - 2$.  Furthermore, on $B_2^c$, we have $W(0, (\ell - 1)r - 1) \cap A_{\ell-1}' = W(0, T) \cap A_{\ell-1}'$.  These results, combined with (\ref{Weq}) when $j = \ell - 1$, imply the theorem.

It remains to prove the claim by induction on $j$.  The claim is trivial when $j = 1$.  Suppose the claim holds for $j-1$, where $j-1 \leq \ell - 2$.  Suppose $k \in S_j \cap J_j^c$.  There are three cases:
\begin{itemize}
\item Suppose $k \leq j-2$.  Then $k \in S_{j-1} \cap J_{j-1}^c$, so by the induction hypothesis it suffices to show that $W(0, (j-1)r - 1) \cap A_k' = W(0, jr - 1) \cap A_k'$.  We proceed by contradiction.  Suppose this equality fails.  Then $X_t = X_s$ for some $t \in A_j$ and $s \in W(0, (j-1)r - 1) \cap [0, kr-w-1]$.  On $B_4^c$, we must also have $s \in A_i'$ for some $i = 1, \dots, k$, and on $B_3^c$, we can not have $i \in J_{j-1}$.  By the induction hypothesis $W(0, (j-1)r - 1) \cap A_i' = \emptyset$ when $i \in S_{j-1}^c \cap J_{j-1}^c$, so $i \in S_{j-1} \cap J_{j-1}^c$.  Using the induction hypothesis again, we get that $W(0, (j-1)r - 1) \cap A_i' = W((i-1)r, ir - 1) \cap A_i'$, which implies that $X_s \in LE(X[A_i])$.  However, this means that $I_{i,j} = 1$, which contradicts that $k \in S_j$.

\item Suppose $k = j-1$.  Again we must show $W(0, (j-1)r - 1) \cap A_k' = W(0, jr - 1) \cap A_k'$.  If the equality fails, then $X_t = X_s$ for some $t \in A_j$ and $s \in W(0, (j-1)r - 1) \cap [0, (j-1)r-w-1]$.  On $B_4^c$, we must either have $s \in A_i'$ for some $i = 1, \dots, j-2$, which leads to a contradiction as in the previous case, or $(j-2)r \leq s \leq (j-1)r - w - 1$, which can not happen on $B_2^c$.

\item Suppose $k = j$.  On $B_1^c$, between times $(j-1)r$ and $(j-1)r + w$, the walk $(X_t)_{t=0}^T$ has a local cutpoint, which we call $u$.  We claim that
\begin{equation}
W(0, jr-1) \cap \{u, \dots, jr-1\} = W((j-1)r, jr-1) \cap \{u, \dots, jr-1\},
\label{Wclaim}
\end{equation}
which will imply (\ref{Weq}).
We prove (\ref{Wclaim}) by contradiction.  Suppose (\ref{Wclaim}) fails.  Then there exist $s$ and $t$ such that $X_t = X_s$, where $u < t \leq jr-1$ and $s \in W(0, u-1) \cup W((j-1)r, u-1)$.  However, since $u$ is a local cutpoint, this can only happen if $t \geq s + 2 \tau_n$, and then on $B_2^c$, it must be that $s \leq (j-2)r - 1$.  However, on $B_4^c$, we must have $s \in A_i'$ for some $i = 1, \dots, j-2$.  As in the first case, this would imply that $I_{i,j} = 1$, contradicting that $k \in J_k^c = J_j^c$.
\end{itemize}

Next, suppose $k \in S_j^c \cap J_j^c$.  Then $I_{h,i} = 1$ for some $h$ and $i$ such that $h < k < i \leq j$ and $h \in S_{i-1}$.  This means that $X_t \in LE(X[A_h])$ for some $t \in A_i$.  Choose $t = \min\{u \in A_i: X_u \in LE(X[A_h])\}$.  On $B_4^c$, it follows that $X_t = X_s$ for some $s \in W((h-1)r, hr - 1) \cap A_h'$.  On $B_3^c$, we must have $h \in J_{i-1}^c$, so by the induction hypothesis, $s \in W(0, (i-1)r - 1) \cap A_h'$.  It follows that $W(0, t) \cap (s, t) = \emptyset$ and since $jr - 1 \geq t$, we have $W(0, jr-1) \cap (s,t) = \emptyset$.  Also, $A_k \subset (s,t)$, so $W(0, jr - 1) \cap A_k = \emptyset$.  This completes the proof of the claim.
\end{proof}

For $i = 1, \dots, \ell$, let $N_j$ be the cardinality of the set $W(0, T) \cap A_j$, which is the number of points from the $j$th segment retained after loop erasure. 
The next proposition bounds the random variables $N_j$.

\begin{Prop}
Suppose $X[0,T]$ is good.  Then the following hold:
\begin{enumerate}
\item If $j \in S \cap J^c$, then $|LE(X[A_j])| - 2w \leq N_j \leq |LE(X[A_j])| + 2w$.

\item If $j \in S^c \cap J^c$, then $N_j = 0$.

\item If $j \in J$, then $0 \leq N_j \leq |LE(X[A_j])| + 2w$.
\end{enumerate}
\label{pathprop}
\end{Prop}

\begin{proof}
Parts 1 and 2 of the proposition follow immediately from Lemma \ref{goodlem}.
It remains to prove part 3.  Suppose $j \in J$.  On $B_3^c$, there can be at most one value of $k$ such that $I_{j,k}$ or $I_{k,j}$ is nonzero, so there are three cases:

\begin{itemize}
\item Suppose $j = \ell$.  On $B_1^c$, the random walk $X[0,T]$ has a local cutpoint at some time between $(\ell-1)r$ and $(\ell-1)r + w$, which we call $u$.  Then
\begin{equation}
W(0, T) \cap \{u, \dots, T\} = W((\ell-1)r, T) \cap \{u, \dots, T\}
\label{newW}
\end{equation}
unless there exist two times $s$ and $t$ such that $X_t = X_s$, where $u < t \leq T$ and $s \in W(0, u-1) \cup W((\ell-1)r, u-1)$.  However, since $u$ is a local cutpoint, this would imply $t \geq s + 2 \tau_n$, which is impossible on $B_2^c \cap B_5^c$.  Therefore, (\ref{newW}) holds.   It follows that $N_{\ell} \leq |LE(X[(\ell-1)r, T])| + w$.  Now if $|LE(X[(\ell - 1)r, T])| > |LE(X[A_{\ell}])| + 2 \tau_n$, then there exist $s < T - 2 \tau_n$ and $t > T$ such that $X_s = X_t$, which is impossible on $B_2^c$.  It follows that $N_{\ell} \leq |LE(X[A_{\ell}])| + w + 2 \tau_n \leq |LE(X[A_{\ell}])| + 2w$.

\item Suppose $I_{k,j} = 1$ for some $k < j < \ell$.  Let $$t = \min\{u \geq (j-1)r: X_v \neq X_s \mbox{ for all }v \in [u, jr - 1] \mbox{ and }s \in [0, kr - 1] \cap W(0, u-1)\}$$ be the first time within the $j$th segment after which there are no more intersections involving earlier segments.  Then $W(0, jr-1) \cap A_j = W(t - 1, jr-1)$.  On $B_1^c$, the random walk $X[0,T]$ has a local cutpoint $u$ in $[t - 1, t - 1 + w]$.  On $B_2^c$, we have $W((j-1)r, jr-1) \cap \{u, \dots, jr-1\} = W(t-1, jr-1) \cap \{u, \dots, jr-1\}$.  Since $W(0, T) \cap A_j \subset W(0, jr-1) \cap A_j$, it follows that $N_j \leq |LE((X_t)_{t \in A_j})| + w$.

\item Suppose $I_{j,k} = 1$ for some $j < k < \ell$.  On $B_3^c$, we must have $j \in J_{k-1}^c$.  By the claim in the proof of Lemma \ref{goodlem}, either $W(0, (k-1)r - 1) \cap A_j = \emptyset$ or $W(0, (k-1)r - 1) \cap A_j' = W((j-1)r, jr-1) \cap A_j'$.  It follows that $W(0, T) \cap A_j' \subset W((j-1)r, jr-1) \cap A_j'$, and therefore $N_j \leq |LE((X_t)_{t \in A_j})| + 2w$.

\end{itemize}
Because these three cases cover all $j \in J$, the proof is complete.
\end{proof}

\begin{Rmk}
{\em To prepare for the construction of the spanning tree in the next section, note that we can create a path whose vertex set is a subset of $\{1, \dots, \ell\}$ in two ways.  First, let ${\cal P}_0 = \{0\}$.  Given the path ${\cal P}_{j-1}$ for some $j \in \{1, \dots, \ell\}$, define ${\cal P}_j$ by the following rules:
\begin{itemize}
\item If there is a $k$ in the path ${\cal P}_{j-1}$ such that $I_{k,j} = 1$ (note that there can be at most one such $k$ on $B_3^c$), then remove the vertices $k+1, \dots, j-1$ from the path, and give the vertex $k$ the new label $j$.

\item Otherwise, add the vertex $j$ to the graph with an edge between $j-1$ and $j$.
\end{itemize}
Let ${\cal P} = {\cal P}_{\ell}$.  Then, by the construction, the vertices in ${\cal P}$ are precisely the points in ${\cal S}_{\ell}$.  If $i_1 < i_2 < \dots < i_k$ are the vertices, ranked in decreasing order, then there is an edge from $i_{j-1}$ to $i_j$ for $j = 1, \dots, k-1$.

Alternatively, consider $LE(X[0,T])$, and give the points in $(X_t)_{t \in W(0,T) \cap A_j}$ the label $j$.  Let ${\cal P}'$ be the path whose vertex set consists of all $j$ such that $N_j > 0$, and such that $i$ and $j$ are connected by an edge if and only if, when we write the points in $W(0,T)$ as $\sigma_1^{0,T} < \sigma_2^{0,T} < \dots < \sigma_{\ell}^{0,T}$, we have $\sigma_k^{0,T} \in A_i$ and $\sigma_{k+1}^{0,T} \in A_j$ for some $k$.  That is, $i$ and $j$ are connected by an edge if there are adjacent points in $LE(X[0,T])$ with one labeled $i$ and one labeled $j$.  By Proposition \ref{pathprop}, the vertex $j$ is in ${\cal P}'$ if $j \in S \cap J^c$ but not if $j \in S^c \cap J^c$.  Therefore, the only difference between the paths ${\cal P}$ and ${\cal P}'$ is that some $j \in J$ could be included in one path but not the other.}
\label{mainrmk}
\end{Rmk}

\section{Coupling with the complete graph}

In this section, we obtain the scaling limit of the uniform spanning tree on $\Z^4_n$ by using Wilson's algorithm and coupling the construction of the uniform spanning tree on $\Z^4_n$ with the construction of the uniform spanning tree on the complete graph.  This is essentially the same as the approach in \cite{perrev}, although we describe the coupling slightly differently by coupling two labeled trees.

Recall that $$r = \lfloor n^2 (\log n)^{9/22} \rfloor,$$ which is the length of the random walk segments on $\Z^4_n$ that we will consider.  Also, let $0 < \theta < 1/22$ be a small positive number, and let $$w = \lfloor n^2 (\log n)^{\theta} \rfloor.$$  The significance of $w$ is that a random walk segment of length at least $w$ has a local cutpoint with high probability (see Proposition \ref{lcprop}).  When necessary, we will assume $n$ is large enough that $w \geq 2 \tau_n$.

The constructions of the spanning trees on the complete graph and on $\Z_n^4$ will require introducing a considerable amount of notation.  For the convenience of the reader, we summarize some of this notation in the table below.  More precise definitions of these quantities will appear when they are introduced.

\bigskip
\begin{tabular}{ll}
$K_{m, \alpha}$ & Complete graph on $m$ vertices with root added \\
$G_{n, \beta}$ & Graph $\Z^4_n$ with root added \\
$({\tilde T}, {\tilde d})$ & Uniform spanning tree on complete graph $K_m$ \\
$(\tilde{\cal{T}}^*, {\tilde d}^*)$ & Labeled spanning tree on $K_{m, \alpha}$ \\
$({\cal T}, d)$ & Uniform spanning tree on $\Z^4_n$ \\
$({\cal T}_k, d')$ & Tree on $\Z^4_n$ built up after $k$ walks in Wilson's algorithm \\
$({\cal{T}}^*, d^*)$ & Labeled spanning tree on $G_{n, \beta}$, after collapsing segments \\
$R_i$ & The $i$th segment of the tree considered in the algorithm \\
$R_i'$ & The interior of the $i$th segment considered in the algorithm \\
$N_j$ & The vertices of the spanning tree on the torus labeled $j$ \\
${\tilde S}_j$ & Vertices in the tree on $K_{m, \alpha}$ after $j$ steps of Wilson's algorithm \\
$S_j$ & Segments in the tree on $G_{n, \beta}$ after $j$ steps of Wilson's algorithm
\end{tabular}

\subsection{A partial spanning tree on the complete graph}

To make the coupling work, we will choose the number of vertices in the complete graph to be
\begin{equation}
m = \lfloor a_n^{-1} (\log n)^{2/11} \rfloor,
\label{mdef}
\end{equation}
where $a_n$ is the constant from Proposition \ref{maincapprop} when we use $r - 2w$ in place of $M$ and $r$ in place of $L$.
Recall from the introduction that $K_m$ is the complete graph on $m$ vertices and $K_{m, \alpha}$ is obtained from $K_m$ by adding a root vertex $\rho$ that is connected to every other vertex by an edge of weight $m/(\alpha \sqrt{m} - 1)$.

We now define a portion of a spanning tree on $K_{m, \alpha}$ using Wilson's algorithm.
Choose vertices $y_1, \dots, y_k$ uniformly at random from $K_m$.  Define trees $\tilde{\cal T}_0, \tilde{\cal T}_1, \dots, \tilde{\cal T}_k$ by running Wilson's algorithm on $K_{m, \alpha}$ starting from the points $\rho, y_1, \dots, y_k$.  More precisely, for $i = 1, \dots, k$, let $Y^i = (Y^i_t)_{t=0}^{\infty}$ be a weighted random walk on $K_{m, \alpha}$ such that $Y^i_0 = y_i$.  Let $\tilde{\cal T}_0$ consist of the single vertex $\rho$.  For $i = 1, \dots, k$, let $V_i = \min\{t: Y^i_t \mbox{ is a vertex of }\tilde{\cal T}_{i-1}\}$ and define $\tilde{\cal T}_i$ by adjoining to $\tilde{\cal T}_{i-1}$ the path $LE(Y^i[0, V_i])$.  Let $\kappa_1 = 1$ and for $1 \leq i \leq k$, let $\kappa_{i+1} = \kappa_i + V_i + 1$.  Let $\ell' = \kappa_{k+1} - 1$.  If $1 \leq j \leq \ell'$, let $f(j)$ be the value of $i$ such that $\kappa_i \leq j < \kappa_{i+1}$, so the $j$th vertex visited in Wilson's algorithm is part of the walk $Y^{f(j)}$.  Then let $v_j = Y^{f(j)}_{j - \kappa_{f(j)}}$, which is the $j$th vertex visited during the construction.

For integers $i$ and $j$ with $1 \leq i < j \leq \ell'$, let $\tilde{I}_{i,j}$ be the indicator of the event that $v_i = v_j \neq \rho$.  For $1 \leq j \leq \ell'$, let $\tilde{I}_{0,j}$ be the indicator of the event that $v_j = \rho$.  We now obtain a tree $\tilde{\cal T}^*$, which will be the same as the tree $\tilde{\cal T}_k$ but with the vertices labeled in the order that they are visited during the construction.  To obtain this tree, we follow Wilson's algorithm one step at a time.  At each step, we keep track of the set ${\tilde S}_j$, which will be the nonzero vertex labels after $j$ steps in the construction.

\bigskip
\noindent {\bf Construction}:  Given $\kappa_1 < \kappa_2 < \dots < \kappa_{k+1}$ and given indicator random variables ${\tilde I}_{i,j}$ for $0 \leq i < j \leq \ell'$, where $\ell' = \kappa_{k+1}-1$, we define a labeled tree $\tilde{\cal T}^*$ as follows.  Begin with just the root vertex, labeled $0$.  Define $\tilde{S}_0 = \emptyset$.  For $j = 1, 2, \dots, \ell'$, let $f(j)$, as above, be the value of $i$ such that $\kappa_i \leq j < \kappa_i$, and proceed inductively as follows.
\begin{itemize}
\item If ${\tilde I}_{i,j} = 0$ for all $i \in \tilde{S}_{j-1}$, then add the vertex $j$ to the graph, and draw an edge from $j-1$ to $j$ unless $j = \kappa_h$ for some $h$.  Also, let $\tilde{S}_j = \tilde{S}_{j-1} \cup \{j\}$.

\item Otherwise, let $i = \min\{h \in \tilde{S}_{j-1}: \tilde{I}_{h,j} = 1\}$.  There are three cases:
\begin{enumerate}
\item If $i = 0$, then draw an edge from $j-1$ to $0$ and let $\tilde{S}_j = \tilde{S}_{j-1}$.

\item If $f(i) = f(j)$, then a loop has formed, so we change the label of the vertex labeled $i$ to $j$ and erase vertices $i+1, \dots, j-1$.  Let $\tilde{S}_j = (\tilde{S}_{j-1} \cap \{1, \dots, i-1\}) \cup\{j\}$.

\item If $f(i) < f(j)$, then the walk has hit the previous tree, so we draw an edge from $j-1$ to $i$ unless $j = \kappa_h$ for some $h$.  Let $\tilde{S}_{j} = \tilde{S}_{j-1}$.
\end{enumerate}
\end{itemize}
The tree $\tilde{\cal T}^*$ is the tree with vertex set $\tilde{S}_{\ell'}$ that has been obtained after $\ell$ steps of this algorithm.

\bigskip
Because $\tilde{\cal T}^*$ is obtained by Wilson's algorithm, it is the same tree as $\tilde{\cal T}_k$ with the vertices relabeled.  Note that the vertex $y_i$ initially has label $\kappa_i$, although the label could change if this vertex is visited another time.  For each $i$ and $j$, there is a unique path in $\tilde{\cal T}^*$ from $y_i$ to $y_j$.  Define $\tilde{d}^*(y_i, y_j)$ to be the number of vertices along this path, unless the path goes through $\rho$ (i.e. the vertex labeled $0$), in which case let $\tilde{d}^*(y_i, y_j) = \infty$.

\subsection{A partial spanning tree on the torus}

We now define a portion of a spanning tree on $G_{n, \beta}$, which is the graph obtained by adding to $\Z^4_n$ a root vertex $\rho$, connected to all other vertices by an edge of weight $8/(\beta n^2 (\log n)^{1/2} - 1)$.  Choose points $x_1, \dots, x_k$ uniformly at random from $\Z^4_n$.  Define trees ${\cal T}_0, {\cal T}_1, \dots, {\cal T}_k$ by running Wilson's algorithm on $G_{n, \beta}$, starting from the points $\rho, x_1, \dots, x_k$.  More precisely, 
for $i = 1, \dots, k$, let $X^i = (X^i_t)_{t=0}^{\infty}$ be a weighted random walk on $G_{n, \beta}$ such that $X^i_0 = x_i$.  We will also adopt the convention that $\rho$ is an absorbing state, so if $X^i_t = \rho$ then we set $X^t_s = \rho$ for all $s \geq t$.  Let ${\cal T}_0$ consist of the single vertex $\rho$.  For $1 \leq i \leq k$, let $U_i = \min\{t: X^i_t \mbox{ is a vertex of }{\cal T}_{i-1}\}$.  Then define ${\cal T}_i$ by adjoining to ${\cal T}_{i-1}$ the path ${\cal P}_i = LE(X^i[0, U_i])$.

To break up the branches of the tree into segments of length $r$, let $\zeta_1 = 1$ and, for $1 \leq i \leq k$, let $\zeta_{i+1} = \zeta_i + \lfloor U_i/r \rfloor + 1$.  Let $\ell = \zeta_{k+1} - 1$.  If $1 \leq j \leq \ell$, let $g(j)$ be the value of $i$ such that $\zeta_i \leq j < \zeta_{i+1}$.  This means that the $j$th segment of length $r$ is part of the walk $X^{g(j)}$.  If $g(j) = i$, define $A_j = \{(j - \zeta_i)r, \dots, (j - \zeta_i + 1)r - 1\}$.  Also, let $A_j' = \{(j - \zeta_i)r + w, \dots, (j - \zeta_i + 1)r - w - 1\},$ which is obtained by removing the first $w$ and the last $w$ points from $A_j$.  To simplify notation, for $1 \leq j \leq \ell$, we will denote the path $X^{g(j)}[A_j]$ by $R_j$, and we denote the path $X^{g(j)}[A_j']$ by $R_j'$.  Also, denote $X^{g(j)}[A_j \setminus A_j']$ by $R_j^-$.  We call the paths $R_j$ segments, and we call two segments $R_{j_1}$ and $R_{j_2}$ adjacent if $g(j_1) = g(j_2)$ and $|j_1 - j_2| = 1$.

We now define events $D_1, \dots, D_4$.  When each of the walks $X^i[0, U_i]$ is good and none of the events $D_1, \dots, D_4$ occurs, it will be possible to analyze the spanning tree by breaking up the random walks into segments of length $r$.

\begin{itemize}
\item Let $D_1$ be the event that some segment intersects two segments not adjacent to it at points other than $\rho$.

\item Let $D_2$ be the event that there are non-adjacent segments $R_{j_1}$ and $R_{j_2}$ such that $R_{j_1}$ and $R_{j_2}^-$ intersect at some point other than $\rho$.

\item Let $D_3$ be the event that some segment intersects $\rho$ and also intersects a non-adjacent segment at a point other than $\rho$.

\item Let $D_4$ be the event that for some $i$, $R_{\zeta_i}$ intersects either $\rho$ or a non-adjacent segment.
\end{itemize}

\begin{Prop}
Let $G$ be the event that $X^i[0, U_i]$ is good for $i = 1, \dots, k$ and that none of the events $D_1$, $D_2$, $D_3$, or $D_4$ occurs.  Then there is a constant $C$, which depends on $k$, such that $P(G) \geq 1 - C (\log n)^{-1/11}$.
\label{goodprop}
\end{Prop}

\begin{proof}
Because the stopping times $U_i$ depend on the paths, we will work also with the times $V_i$, where $V_i = \min\{t: X^i_t = \rho\}$.  Then the stopping times $V_1, \dots, V_k$ are independent of the walks $X^1, \dots, X^k$.  Each $V_i$ has a geometric distribution with mean $\beta n^2 (\log n)^{1/2}$, so there is a constant $C$ such that $E[V_i^j] \leq C n^{2j} (\log n)^{j/2}$ for $i = 1, \dots, k$ and $j = 1, 2, 3$.  Also, since $\rho$ is in the tree ${\cal T}_0$, we have $U_i \leq V_i$ for all $i$.

Because we will work with the potentially longer walks $X^i[0, V_i]$ in the proof, it is necessary to relabel the segments.  Let $\bar{\zeta}_1 = 1$ and, for $1 \leq i \leq k$, let $\bar{\zeta}_{i+1} = \bar{\zeta}_i + \lfloor V_i/r \rfloor + 1$.  Let $\bar{\ell} = \bar{\zeta}_{k+1} - 1$.  If $1 \leq j \leq \bar{\ell}$, let $\bar{g}(j)$ be the value of $j$ such that $\bar{\zeta}_i \leq j < \bar{\zeta}_{i+1}$.  If $\bar{g}(j) = i$, define $\bar{A}_j = \{(j - \bar{\zeta}_i)r, \dots, (j - \bar{\zeta}_i + 1)r - 1\}$ and $\bar{A}_j' = \{(j - \bar{\zeta}_i)r + w, \dots, (j - \bar{\zeta}_i + 1)r - w - 1\}$.  For $1 \leq j \leq \bar{\ell}$, define $\bar{R}_j$, $\bar{R}_j'$, and $\bar{R}_j^-$ in the same way as $R_j$, $R_j'$, and $R_j^-$ but using $\bar{g}$, $\bar{\zeta}_i$, $\bar{A}_j$, and $\bar{A}_j'$ in place of $g$, $\zeta_i$, $A_j$, and $A_j'$ respectively.  Finally, for $i = 1, 2, 3, 4$, let $\bar{D}_i$ be defined in the same way as $D_i$.  Because $U_i \leq V_i$ for all $i$, we have $D_i \subset \bar{D}_i$ for $i = 1, 2, 3, 4$, so to bound $P(D_i)$ it will suffice to bound $P(\bar{D}_i)$.

Recall that five conditions must be satisfied for the walk $X^i[0, U_i]$ to be good.  Define the events $B_1^i, \dots, B_5^i$ in the same way as the events $B_1, \dots, B_5$ were defined in the previous section, but using $X^i[0, U_i]$ in place of $X[0, T]$.
We first consider $B_1^i$.  By Proposition \ref{lcprop}, for any $j \in \N$, the probability that a random walk $(X_t)_{t=0}^{\infty}$ on $\Z^4_n$ contains an interval of the form
$[t, t + w]$ with $(j-1)n^2 (\log n)^{1/2} \leq t < j n^2 (\log n)^{1/2}$ that has no local cutpoint is at most $C/(\log n)$.  Since $U_i \leq V_i$ and $V_i$ is independent of $X^i$, it follows that
$$P(B_1^i|V_i) \leq \frac{C}{\log n} \bigg( \bigg\lfloor \frac{V_i}{n^2 (\log n)^{1/2}} \bigg\rfloor + 1 \bigg).$$  By taking expectations of both sides, we get $P(B_1^i) \leq C/(\log n)$ for $i = 1, \dots, k$.

We next consider the event $B_2^i$.  If $\zeta_i \leq j < \zeta_{i+1}$, then the walk $X^i[A_j \cup A_{j+1}]$ has length $2r$.  If $(X_t)_{t=0}^{2r}$ is a random walk on $\Z^4_n$, then by Proposition \ref{medloops}, the probability that there are two times $s$ and $t$ such that $|s - t| \geq 2 \tau_n$ and $X_s = X_t$ is at most $C (2r)^2/(n^4 \log n)$.  It follows that $$P(B_2^i|V_i) \leq \bigg\lfloor \frac{V_i}{r} \bigg\rfloor \frac{4C r^2}{n^4 \log n} \leq \frac{4C V_i r}{n^4 \log n}.$$
Taking expectations of both sides gives $P(B_2^i) \leq C (\log n)^{-1/11}$.

To show that each of the walks $X^i[0, U_i]$ is good with high probability, it remains to bound $P(B_3^i)$, $P(B_4^i)$, and $P(B_5^i)$.  However, if $B_3^i$ occurs for some $i$, then $D_1$ occurs, and if $B_4^i$ occurs for some $i$, then $D_2$ occurs.  Also, the last segment of $X^i[0, U_i]$ must either hit $\rho$ or one of the segments $X^j[0, U_j]$ for $j < i$.  Therefore, on $D_1^c \cap D_3^c$, the event $B_5^i$ can not occur.  Thus, to complete the proof, it suffices to show that $P(\bar{D}_i) \leq C(\log n)^{-1/11}$ for $i = 1, 2, 3, 4$.

We first bound $P(\bar{D}_2)$.  Let $W_{j_1, j_2}$ be the event that $\bar{R}_{j_1}$ and $\bar{R}_{j_2}^-$ intersect at some point other than $\rho$, meaning that $\bar{R}_{j_1}$ intersects either the first $w$ or the last $w$ steps of $\bar{R}_{j_2}$.  The segments $\bar{R}_{j_1}$ and $\bar{R}_{j_2}$ are either parts of different walks or are separated in time by at least $\tau_n$, so Lemma \ref{mixlem} and Proposition \ref{mainintprop} imply that $P(W_{j_1, j_2}) \leq Crw/(n^4 \log n)$.
There are at most $\bar{\ell}^2$ pairs $(j_1, j_2)$ to be considered, so $$P(\bar{D}_2|V_1, \dots, V_k) \leq \frac{Crw \bar{\ell}^2}{n^4 \log n}.$$ By taking expectations of both sides, and using that $E[\bar{\ell}^2] \leq C n^4 r^{-2} \log n$ for some constant $C$, we see that there is a constant $C$ depending on $k$ such that $P(\bar{D}_2) \leq C (\log n)^{-9/22 + \theta}$.

To bound $P(\bar{D}_1)$, suppose $j_1, j_2, j_3 \leq \bar{\ell}$ are distinct integers and neither $R_{j_2}$ nor $R_{j_3}$ is adjacent to $R_{j_1}$.  
On the event $\bar{D}_2^c$, the walk $\bar{R}_{j_1}$ can only intersect $\bar{R}_{j_2}$ and $\bar{R}_{j_3}$ if the walk $\bar{R}_{j_1}'$ intersects both $\bar{R}_{j_2}'$ and $\bar{R}_{j_3}'$, an event that we call $W_{j_1, j_2, j_3}$.  Let $\lambda_1$, $\lambda_2$, and $\lambda_3$ be any paths in $\Z^4_n$ of length $r - 2w$, and let $S_1$, $S_2$, and $S_3$ be paths in $\Z^4_n$ that come from independent random walks of length $r - 2w$, started from the uniform distribution.  Since the segments $\bar{R}_{j_1}'$, $\bar{R}_{j_2}'$, and $\bar{R}_{j_3}'$ are separated from one another by at least $2w \geq \tau_n$, applying Lemma \ref{mixlem} twice gives
\begin{equation}
P(\bar{R}_{j_1}' = \lambda_1, \bar{R}_{j_2}' = \lambda_2, \bar{R}_{j_3}' = \lambda_3) \leq \frac{9}{4} P(S_1 = \lambda_1, S_2 = \lambda_2, S_3 = \lambda_3).
\label{Rlambda}
\end{equation}
Therefore, the probability that $\bar{R}_{j_1}'$ intersects both $\bar{R}_{j_2}'$ and $\bar{R}_{j_3}'$ is at most $9/4$ times the probability that $S_1$ intersects $S_2$ and $S_3$.  Thus, by Proposition \ref{closeprop},
\begin{equation}
\label{Wj123}
P(W_{j_1, j_2, j_3}) \leq \frac{9}{4} E[\mbox{Close}_{r-2w}(S_2, S_3)] \leq \frac{C r^4}{n^8 (\log n)^2}.
\end{equation}
Since there are at most $\bar{\ell}^3$ possible choices of $(j_1, j_2, j_3)$, we get $$P(\bar{D}_1|V_1, \dots, V_k) \leq \frac{C r^4 \bar{\ell}^3}{n^8 (\log n)^2}.$$  Since $E[\bar{\ell}^3] \leq C n^6 r^{-3} (\log n)^{3/2}$ for some constant $C$, taking expectations of both sides gives $P(\bar{D}_1) \leq C (\log n)^{-1/11}$.

For $\bar{D}_3$, note that if $\bar{R}_{j_1}$ intersects $\rho$, then $j_1 = \bar{\zeta}_{i+1}-1$ for some $i = 1, \dots, k$.  If $R_{j_1}$ and $R_{j_2}$ are not adjacent segments, then by Proposition \ref{mainintprop}, the probability that $\bar{R}_{j_1}$ and $\bar{R}_{j_2}$ intersect at some point other than $\rho$ is at most $C r^2/(n^4 \log n)$.  Conditional on $V_1, \dots, V_k$, there are $k$ possible choices for $j_1$ and at most $\bar{\ell}$ for $j_2$, so $$P(\bar{D}_3|V_1, \dots, V_k) \leq \frac{C r^2 \bar{\ell}}{n^4 \log n}.$$  Since $E[\bar{\ell}] \leq C n^2 r^{-1} (\log n)^{1/2}$, we take expectations of both sides to get $P(\bar{D}_3) \leq C (\log n)^{-1/11}$.

Finally, to bound $P(\bar{D}_4)$, suppose $j_1 = \bar{\zeta}_i$.  Because the time that the random walk $X^i$ first visits $\rho$ has a geometric distribution with mean $\beta n^2 (\log n)^{1/2}$, the probability that $\bar{R}_{j_1}$ intersects $\rho$ is at most $r/(\beta n^2 (\log n)^{1/2}) \leq C(\log n)^{-1/11}$.  If $j_2 \notin \{j_1, j_1 + 1\}$, then the probability that the walks $\bar{R}_{j_1}$ and $\bar{R}_{j_2}$ intersect is at most $C r^2/(n^4 \log n)$ by Proposition \ref{mainintprop}.  Therefore, following the same reasoning used to bound $P(\bar{D}_3)$, we get $P(\bar{D}_4) \leq C (\log n)^{-1/11}$.
\end{proof}

For the rest of this subsection, we will work on the event $G$.  The next step is to label the vertices of the tree ${\cal T}_k$ by considering the segments of length $r$ in the order that they are visited by Wilson's algorithm.  To define the time indices that are retained after loop-erasure of the walk $X^i[0, U_i]$, let $\sigma^{u,v,i}_0 = \max\{t: u \leq t \leq v \mbox{ and }X_t = X_u\}$.  For $j \leq 1$, if $\sigma^{u,v,i}_{j-1} < v$, then let $$\sigma^{u,v,i}_j = \max\{t: \sigma_{j-1}^{u,v,i} < t \leq v \mbox{ and }X_t = X_{\sigma_{i-1}^{u,v,i} + 1}\}.$$  Let $W^i(u, v)$ denote the set of all times $\sigma^{u,v,i}_j$, which are the times retained in $LE(X^i[u, v])$.  Give the root the label $0$.  For $1 \leq i \leq k$ and $\zeta_i \leq j < \zeta_{i+1}$, give the label $j$ to all vertices of the tree ${\cal T}_k$ that equal $X^i_t$ for some $t$ in $W^i(0, U_i) \cap A_j$, other than the vertex $X^i_{U_i}$.  Informally, the label $j$ is given to all vertices in ${\cal T}_k$ that are in the $j$th segment of length $r$ visited by Wilson's algorithm.  For $j = 1, \dots, \ell$, let $N_j$ be the number of vertices labeled $j$.

For integers $i$ and $j$ such that $1 \leq i < j \leq \ell$ and $R_i$ and $R_j$ are not adjacent, define $I_{i,j}$ to be the indicator of the event that $LE(R_i) \cap R_j' \neq \emptyset$.  Note that we use $R_j'$ for the definition rather than $R_j$.  This will be convenient for the proof of Proposition \ref{coupleprop} later, but on $D_2^c$ the two definitions are equivalent.  Also, for $1 \leq j \leq \ell$, let $I_{0,j}$ be the indicator of the event that $R_j$ intersects $\rho$, and if $j \geq 2$ and $g(j) = g(j-1)$, let $I_{j-1, j} = 0$.  Note that on $D_1^c \cap D_3^c$, for $1 \leq j \leq \ell$, at most one of the random variables $I_{0,j}, I_{1,j}, \dots, I_{j-1, j}$ can be $1$.

We can use these indicator random variables to construct a tree ${\cal T}^*$.  Define ${\cal T}^*$ to be the tree obtained from the construction in subsection 4.1, with $\zeta_1, \dots, \zeta_{k+1}$ in place of $\kappa_1, \dots, \kappa_{k+1}$ and the $I_{i,j}$ in place of ${\tilde I}_{i,j}$.  Also, denote by $S_0, \dots, S_{\ell}$ the sets that were called ${\tilde S}_0, \dots, {\tilde S}_{\ell'}$ in subsection 4.1, and let $S = S_{\ell}$.  We can view ${\cal T}^*$ as being the tree obtained by collapsing each segment of length $r$ in ${\cal T}_k$ down to a single vertex.  The vertices of ${\cal T}^*$ correspond to segments of length $r$ that are retained after loop erasure, with the convention that when segment $j$ intersects segment $i$ with $i < j$, segment $j$ is retained while segment $i$ is not.
Also, $S_j$ consists of segments that have not been erased after $j$ segments have been processed by the algorithm.  Finally, let $J_j = \{i \leq j: I_{i,h} = 1 \mbox{ or }I_{h,i} = 1 \mbox{ for some }0 \leq h \leq j\}$, and let $J = J_{\ell}$.  Thus, $J$ consists of the segments that are involved in intersections.

The next lemma shows that in the construction of ${\cal T}^*$, for $i = 1, \dots, k$, the vertex $\zeta_{i+1}-1$ intersects the tree that has been built up after stage $\zeta_i - 1$, which implies that indeed the final graph ${\cal T}^*$ is a tree.

\begin{Lemma}
On $G$, we have $\zeta_{i+1} - 1 = \min\{j \geq \zeta_i: I_{h,j} = 1 \mbox{ for some }h \in S_{\zeta_i-1} \cup \{0\}\}$ for $i = 1, \dots, k$.
\label{endpointlem}
\end{Lemma}

\begin{proof}
We proceed by induction.  First consider $i = 1$.  We have $U_1 = \min\{t: X^1_t = \rho\}$.  Also, note that $\zeta_1 = 1$ and $S_0 = \emptyset$, and $\min\{j: I_{0,j} = 1\} = 1 + \lfloor U_1/r \rfloor$.  Since, by definition, $\zeta_2 - 1 = 1 + \lfloor U_1/r \rfloor$, the result for $i = 1$ follows.

Now, consider $i > 1$.  The tree ${\cal T}_{i-1}$ is simply a union of the paths ${\cal P}_1, \dots, {\cal P}_{i-1}$, all of which come from loop-erasing good walks.  By Lemma \ref{goodlem}, for $h \leq i - 1$, the only differences between ${\cal P}_h$ and the corresponding loop-erased segments in $S \cap [\zeta_h, \zeta_{h+1} - 1]$ are the segments in $J$ that are involved in intersections, and the first and last $w$ steps of the segments.  That is, for $h = 1, \dots, i-1$ we have $${\cal P}_h \: \triangle \bigg( \bigcup_{j \in S \cap [\zeta_h, \zeta_{h+1}-1]} LE(R_j) \bigg) \subset \bigg( \bigcup_{j \in J_{\zeta_{h+1}-1} \cap[\zeta_h, \zeta_{h+1}-1]} R_j \bigg) \cup \bigg( \bigcup_{j=\zeta_h}^{\zeta_{h+1}-1} R_j^- \bigg).$$
Note that by the strong induction hypothesis, $\zeta_{h+1} - 1 \in J_{\zeta_{h+1}-1}$ for $h = 1, \dots, i-1$, which is important because the last segment of the path is always included in $J$ in the definition used in Lemma \ref{goodlem}.
Because indices from previous paths are never removed from $S$, we have $S \cap [1, \zeta_i - 1] = S_{\zeta_i - 1}$.  It follows that
\begin{equation}
{\cal T}_{i-1} \: \triangle \bigg( \bigcup_{j \in S_{\zeta_i - 1}} LE(R_j) \bigg) \subset \bigg( \bigcup_{j \in J_{\zeta_i-1}} R_j \bigg) \cup \bigg( \bigcup_{j=1}^{\zeta_i-1} R_j^- \bigg).
\label{Tdelta}
\end{equation}

Recall that $U_i = \min\{t: X^i_t \mbox{ is a vertex of }{\cal T}_{i-1}\}$.  On $D_2^c$, we can not have $X^i_{U_i} \in R_j^-$ for any $j < \zeta_i$.  If $j \in J_{\zeta_i - 1}$, then either $R_j$ intersects some non-adjacent segment $R_h$ for some $h \leq \zeta_i - 1$, or $R_j$ intersects $\rho$.  Therefore on $D_1^c \cap D_3^c$, if $j \in J_{\zeta_i - 1}$, the point $X^i_{U_i}$ can not be on the path $R_j'$.  It follows from (\ref{Tdelta}) that either $X^i_{U_i} = \rho$, in which case $I_{0, \zeta_{i+1}-1} = 1$, or $X^i_{U_i}$ is in $LE(R_j)$ for some $j \in S_{\zeta_i-1}$, which on $D_2^c$ means $I_{j,\zeta_{i+1}-1} = 1$.  Therefore, $\zeta_{i+1} - 1 \geq \min\{j \geq \zeta_i: I_{h,j} = 1 \mbox{ for some }h \in S_{\zeta_i-1} \cup \{0\}\}$.

Now, suppose $j \geq \zeta_i$, $h \in S_{\zeta_i - 1} \cup \{0\}$, and $I_{h,j} = 1$.  If $h = 0$, then $R_j$ intersects $\rho$.  If $h \in S_{\zeta_i - 1}$, then $R_j'$ intersects $LE(R_h)$.  It follows from (\ref{Tdelta}) that on $D_1^c \cap D_2^c \cap D_3^c$, the segment $R_j$ intersects ${\cal T}_{i-1}$.  Therefore, $\zeta_{i+1} - 1 \leq h$, which completes the proof.
\end{proof}

\begin{Lemma}
On $G$, the following hold:
\begin{enumerate}
\item If $j \in S \cap J^c$, then $|LE(R_j)| - 2w \leq N_j \leq |LE(R_j)| + 2w$.

\item If $j \in S^c \cap J^c$, then $N_j = 0$.

\item If $j \in J$, then $0 \leq N_j \leq |LE(R_j)| + 2w$.
\end{enumerate}
\label{labelprop}
\end{Lemma}

\begin{proof}
On $G$, the tree ${\cal T}_k$ is a union of the paths ${\cal P}_1, \dots, {\cal P}_k$, all of which come from good walks.  Therefore, this result is an immediate consequence of Proposition \ref{pathprop}, once it is clear that the sets $S$ and $J$ are defined for the tree ${\cal T}_k$ the same manner as they are defined in Proposition \ref{pathprop}.  However, this can be seen from the definition, and the fact that $\zeta_{i+1} - 1 \in J$ for all $i = 1, \dots, k$ by Lemma \ref{endpointlem}.  Additional indices are placed in $J$ because of intersections between different paths, but this does not affect the result because the claims for $j \in J$ are weaker than those for $j \in J^c$.
\end{proof}

If $1 \leq i < j \leq k$, then there is a unique path in ${\cal T}_k$ between the vertices $x_i$ and $x_j$.  Let $d'(x_i, x_j)$ be the number of vertices in this path, unless the path goes through the root $\rho$, in which case let $d'(x_i, x_j) = \infty$.  Note that for $i = 1, \dots, k$, the vertex $x_i$ initially gets the label $\zeta_i$, and on $D_4^c$ this label never changes.  Therefore, the path in ${\cal T}^*$ from $\zeta_i$ to $\zeta_j$ corresponds to the path in ${\cal T}_k$ from $x_i$ to $x_j$.  Define $d^*(x_i, x_j)$ to be the number of vertices of the path in the tree ${\cal T}^*$ between $\zeta_i$ and $\zeta_j$, unless this path goes through $0$, in which case let $d^*(x_i, x_j) = \infty$.  The next proposition allows us to compare $d'(x_i, x_j)$ and $d^*(x_i, x_j)$.

\begin{Prop}
Let $S_{i,j}$ be the set of vertices on the path in ${\cal T}^*$ from $\zeta_i$ to $\zeta_j$.  On $G$, we have $$\bigg| d'(x_i, x_j) - \sum_{h \in S_{i,j}} N_h \bigg| \leq \sum_{h \in J} N_h.$$
\label{distprop}
\end{Prop}

\begin{proof}
We claim that if $h \in J^c$, then the following are equivalent:
\begin{itemize}
\item There is a vertex labeled $h$ on the path in ${\cal T}_k$ from $x_i$ to $x_j$.

\item Every vertex labeled $h$ is on the path in ${\cal T}_k$ from $x_i$ to $x_j$.

\item We have $h \in S_{i,j}$.
\end{itemize}
This claim immediately implies the proposition.  

To prove the claim, we first show that the first two statements are equivalent.  The path in ${\cal T}_k$ from $x_i$ to $x_j$ is made up of portions of the paths ${\cal P}_1, \dots, {\cal P}_k$.  The vertices labeled $h$, if there are any, are consecutive points along the path ${\cal P}_{g(h)}$.  Therefore, the only way that the path in ${\cal T}_k$ from $x_i$ to $x_j$ can traverse some but not all of the vertices labeled $h$ is if one of the other paths ${\cal P}_b$ with $b > g(h)$ hits the path ${\cal P}_{g(h)}$ at one of the vertices labeled $h$.  However, as observed in the proof of Lemma \ref{endpointlem}, on $G$ this implies that $I_{h, \zeta_{b+1}-1} = 1$ and therefore $h \in J$.

It remains to show that the first two statements are equivalent to $h \in S_{i,j}$.  However, this part of the claim is just an extension of Remark \ref{mainrmk}.  Each of the paths ${\cal P}_1, \dots, {\cal P}_k$ comes from loop-erasing a good segment of a random walk, so by Remark \ref{mainrmk}, the only difference between the set of labels of vertices in ${\cal P}_i$ and the set $S \cap [\zeta_i, \zeta_{i+1}-1]$ is that numbers in $J$ could appear in one but not the other.  Because the path in ${\cal T}_k$ from $x_i$ to $x_j$ traverses portions of these paths, and the path in ${\cal T}^*$ from $\zeta_i$ to $\zeta_j$ traverses the corresponding paths in ${\cal T}^*$, it follows that if $h \in J^c$ then $h \in S_{i,j}$ if and only if there is a vertex labeled $h$ on the path in ${\cal T}_k$ from $x_i$ to $x_j$.
\end{proof}

\subsection{Bound on the probability of coupling}

In this section, we show that it is possible to couple the trees $\tilde{\cal T}^*$ and ${\cal T}^*$ such that $\tilde{\cal T}^* = {\cal T}^*$ with high probability.  This will allow us to compare the distances between points in a uniform spanning tree on the complete graph and distances between points in a uniform spanning tree on the torus.  We first state a lemma that we will use to bound the probability that our coupling is sucessful.

\begin{Lemma}
Let $Y_1, \dots, Y_n$ and $Z_1, \dots, Z_n$ be $\{0, 1\}$-valued random variables.  Let $p_i = P(Y_i = 1)$ and $q_i = P(Z_i = 1)$.  Let $p_{ij} = P(Z_i = Z_j = 1)$ and $q_{ij} = P(Y_i = Y_j = 1)$.  Then there exists a coupling of $(Y_1, \dots, Y_n)$ and $(Z_1, \dots, Z_n)$ such that $$P(Y_i \neq Z_i \mbox{ for some }i = 1, \dots, n) \leq \sum_{i=1}^n \sum_{j \neq i} (p_{ij} + q_{ij}) + \sum_{i=1}^n |p_i - q_i|.$$
\label{couplelem}
\end{Lemma}

\begin{proof}
Note that $P(Y_i = 0 \mbox{ for all }i) \geq 1 - \sum_{i=1}^n p_i$ and $P(Z_i = 0 \mbox{ for all }i) \geq 1 - \sum_{i=1}^n q_i$.  For $i = 1, \dots, n$, we have $P(Y_i = 1 \mbox{ and }Y_j = 0 \mbox{ for }j \neq i) \geq p_i - \sum_{j \neq i} p_{ij}$ and $P(Z_i = 1 \mbox{ and }Z_j = 0 \mbox{ for }j \neq i) \geq q_i - \sum_{j \neq i} q_{ij}$.  Therefore, the two sequences $(Y_1, \dots, Y_n)$ and $(Z_1, \dots, Z_n)$ can be coupled so that $Y_i = Z_i = 0$ for all $i$ with probability at least $\min\{1 - \sum_{j=1}^n p_j, 1 - \sum_{j=1}^n q_j\}$ and, for each $i$, we have $Y_i = Z_i = 1$ and $Y_j = Z_j = 0$ for $j \neq i$ with probability at least $\min\{p_i, q_i\} - \sum_{j \neq i} (p_{ij} + q_{ij})$.  With this coupling, if $A$ is the event that $Y_i = Z_i$ for all $i$, then
\begin{align}
P(A) &\geq \min\bigg\{1 - \sum_{i=1}^n p_i, 1 - \sum_{i=1}^n q_i\bigg\} + \sum_{i=1}^n \min\{p_i, q_i\} - \sum_{i=1}^n \sum_{j \neq i} (p_{ij} + q_{ij}) \nonumber \\
&\geq 1 - \sum_{i=1}^n \max\{p_i, q_i\} + \sum_{i=1}^n \min\{p_i, q_i\} - \sum_{i=1}^n \sum_{j \neq i} (p_{ij} + q_{ij}) \nonumber \\
&= 1 - \sum_{i=1}^n |p_i - q_i| - \sum_{i=1}^n \sum_{j \neq i} (p_{ij} + q_{ij}), \nonumber
\end{align}
which implies the lemma.
\end{proof}

\begin{Lemma}
There exist constants $C$ and $C'$, depending on $k$, such that $$P(\ell > C' (\log n)^{2/22} (\log \log n)) \leq \frac{C}{\log n}.$$
\label{Pllem}
\end{Lemma}

\begin{proof}
Because the time for a random walk on $G_{n, \beta}$ to visit the root has a geometric distribution with mean $\beta n^2 (\log n)^{1/2}$, the distribution of $U_1 + \dots + U_k$ is stochastically dominated by the distribution of the sum of $k$ independent random variables, which we call $Z_1, \dots, Z_k$, each having the geometric distribution with mean $\beta n^2 (\log n)^{1/2}$.  If $Z$ has the geometric distribution with mean $m \geq 2$, then
$$P(Z > am) = \bigg(1 - \frac{1}{m} \bigg)^{\lfloor am \rfloor} \leq \bigg(1 - \frac{1}{m} \bigg)^{-1} e^{-a} \leq 2 e^{-a}.$$  Therefore, for sufficiently large $n$, 
\begin{align}
P(U_1 + \dots + U_k > k \beta n^2 (\log n)^{1/2} (\log \log n)) &\leq k P(Z_i > \beta n^2 (\log n)^{1/2} (\log \log n)) \nonumber \\
&\leq 2k e^{-\log \log n} = \frac{2k}{\log n}.
\end{align}
Because $\ell = \sum_{i=1}^k (1 + \lfloor U_i/r \rfloor)$ and $r = \lfloor n^2 (\log n)^{9/22} \rfloor$, the lemma follows.
\end{proof}

\begin{Prop}
For all $\theta > 0$, the random variables $I_{i,j}$ and $\tilde{I}_{i,j}$ can be coupled so that 
$$P(I_{i,j} = \tilde{I}_{i,j} \mbox{ for all }j \leq \ell \mbox{ and }i \in S_{j-1} \cup \{0\}) \geq 1 - \frac{C}{(\log n)^{1/22 - \theta}}$$ for some constant $C$ that depends on $k$.
\label{coupleprop}
\end{Prop}

\begin{Rmk}
{\em Suppose that $I_{i,j} = \tilde{I}_{i,j}$ for all $j \leq \ell$ and $i \in S_{j-1} \cup \{0\}$.  Then, by induction on $j$, it follows that $S_j= \tilde{S}_j$ for all $j \leq \ell$.  Therefore, by the construction of the trees ${\cal T}^*$ and $\tilde{{\cal T}}^*$, we have ${\cal T}^* = \tilde{{\cal T}}^*$.  Thus, Proposition \ref{coupleprop} implies that $P({\cal T}^* = \tilde{{\cal T}}^*) \geq 1 - C (\log n)^{-1/22 + \theta}$.}
\label{Trem}
\end{Rmk}

\begin{proof}[Proof of Proposition \ref{coupleprop}]
Our strategy will be to couple the random variables $I_{i,j}$ and $\tilde{I}_{i,j}$ for fixed $j$, using Lemma \ref{couplelem}, and then proceed by induction on $j$.
We first consider the random variables $\tilde{I}_{i,j}$ conditional on $\tilde{S}_{j-1}$.  Because the random walk on $K_{m, \alpha}$ goes to $\rho$ on each step with probability $1/(\alpha \sqrt{m})$, we have
\begin{equation}
P(\tilde{I}_{0,j} = 1|\tilde{S}_{j-1}) = \frac{1}{\alpha \sqrt{m}}.
\label{tildecond1}
\end{equation}
Also, when the random walk on $K_{m, \alpha}$ does not go to $\rho$, it visits a randomly chosen vertex.  Therefore, if $i \in \tilde{S}_{j-1}$, then
\begin{equation}
P(\tilde{I}_{i,j} = 1|\tilde{S}_{j-1}) = \bigg(1 - \frac{1}{\alpha \sqrt{m}} \bigg) \frac{1}{m}.
\label{tildecond2}
\end{equation}
If $h, i \in \tilde{S}_{j-1}$, then $v_h \neq v_i$.  Because the random walk can only visit one vertex at a time, if $h, i \in \tilde{S}_{j-1}$ and $h \neq i$, then
\begin{equation}
P(\tilde{I}_{h,j} = \tilde{I}_{i,j} = 1|\tilde{S}_{j-1}) = 0.
\label{tildecond3}
\end{equation}

For all $j \geq 0$, let ${\cal F}_j$ be the $\sigma$-field generated by the segments $R_i$ for $i \leq j$.  Note that the event $i \in S_j$ is in ${\cal F}_j$.  Since $I_{0,j} = 1$ if and only if $R_j$ intersects $\rho$, and on each step the random walk $X^{g(j)}$ goes to $\rho$ with probability $1/(\beta n^2 (\log n)^{1/2})$, we have
\begin{equation}
P(I_{0,j} = 1|{\cal F}_{j-1}) = 1 - \bigg(1 - \frac{1}{\beta n^2 (\log n)^{1/2}} \bigg)^r.
\label{condI0j}
\end{equation}
It follows from (\ref{tightmixbound}) that for all $x \in \Z^4_n$, we have
\begin{equation}
\bigg| P(X^{g(j)}_{(j - \zeta_i)r + w} = x|{\cal F}_{j-1}) - \frac{1}{n^4} \bigg| \leq \frac{1}{n^4} 2^{-\lfloor w/\tau_n \rfloor} \leq \frac{C}{n^4 \log n}
\label{startbias}
\end{equation}
for some constant $C$.  If the walk $R_j'$ started from a point that was uniformly distributed on $\Z^4_n$, then since $A_j'$ is a segment of length $r - 2w$, this walk would intersect a set $U$ with probability $\mbox{Cap}_{r - 2w}(U)$ and would intersect both $U$ and $V$ with probability $\mbox{Close}_{r - 2w}(U, V)$.  Combining this observation with (\ref{startbias}) and the definition of $I_{i,j}$, we get that for $i \in S_{j-1}$,
$$P(I_{i,j} = 1|{\cal F}_{j-1}) = \bigg(1 - \frac{1}{\beta n^2 (\log n)^{1/2}} \bigg)^r (1 + E_{i,j}) \mbox{Cap}_{r - 2w}\big(LE(R_i)\big),$$
where the first term is $P(I_{0,j} = 0)$ and the error term $|E_{i,j}|$ is bounded by $C/(\log n)$.  Since $P(I_{0,j} = 1) \leq r/(\beta n^2 (\log n)^{1/2}) \leq C/(\log n)^{1/11}$, we have
\begin{equation}
P(I_{i,j} = 1|{\cal F}_{j-1}) = (1 + E_{i,j}') \mbox{Cap}_{r - 2w}\big(LE(R_i)\big),
\label{condIij}
\end{equation}
where the error term $|E_{i,j}'|$ is bounded by $C/(\log n)^{1/11}$.
Likewise, if $i \in S_{j-1}$, then
\begin{equation}
P(I_{0,j} = I_{i,j} = 1|{\cal F}_{j-1}) \leq \bigg( 1 - \bigg(1 - \frac{1}{\beta n^2 (\log n)^{1/2}} \bigg)^r \bigg) \bigg(1 + \frac{C}{\log n} \bigg) \mbox{Cap}_{r - 2w}\big(LE(R_i)\big)
\label{condI0jIij}
\end{equation}
and if $h, i \in S_{j-1}$ with $h \neq i$, then
\begin{equation}
P(I_{h,j} = I_{i,j} = 1|{\cal F}_{j-1}) \leq \bigg( 1 + \frac{C}{\log n} \bigg) \mbox{Close}_{r - 2w}\big(LE(R_h), LE(R_i)\big).
\label{condIhjIij}
\end{equation}

Let $W_j$ be the event that $$j = \min\{i: I_{h,i} \neq {\tilde I}_{h,i} \mbox{ for some }h \in S_{i-1} \cup \{0\}\}.$$  If none of the events $W_1, \dots, W_{j-1}$ occurs, then $S_{j-1} = \tilde{S}_{j-1}$.  Therefore, by combining equations (\ref{tildecond1}), (\ref{tildecond2}), (\ref{tildecond3}), (\ref{condI0j}), (\ref{condIij}), (\ref{condI0jIij}), and (\ref{condIhjIij}) and applying Lemma \ref{couplelem}, we see that the random variables $I_{i,j}$ and $\tilde{I}_{i,j}$ can be coupled so that 
\begin{align}
P(W_j|{\cal F}_{j-1}) &\leq \bigg| \frac{1}{\alpha \sqrt{m}} - \bigg(1 - \bigg(1 - \frac{1}{\beta n^2 (\log n)^{1/2}} \bigg)^r\bigg) \bigg| \nonumber \\
&\hspace{.1in}+ \sum_{i \in S_{j-1}} \bigg| \bigg(1 - \frac{1}{\alpha \sqrt{m}} \bigg) \frac{1}{m} -
(1 + E_{i,j}') \mbox{Cap}_{r - 2w}\big(LE(R_i)\big) \bigg| \nonumber \\
&\hspace{.1in}+ \bigg(1 + \frac{C}{\log n} \bigg) \bigg( 1 - \bigg(1 - \frac{1}{\beta n^2 (\log n)^{1/2}} \bigg)^r \bigg) \sum_{i \in S_{j-1}} \mbox{Cap}_{r - 2w}\big(LE(R_i)\big)
\nonumber \\
&\hspace{.1in}+ \sum_{h, i \in S_{j-1}, i \neq h} \bigg( 1 + \frac{C}{\log n} \bigg) \mbox{Close}_{r - 2w}\big(LE(R_h), LE(R_i)\big).
\label{bigbound}
\end{align}

By (\ref{alphadef}), the first term on the right-hand side of (\ref{bigbound}) is zero.
To bound the third term, note that $E[\mbox{Cap}_{r - 2w}(LE(R_i))] \leq Cr^2/(n^4 \log n)$ by Proposition \ref{lpssmallprop}, so the expected value of the third term is at most
$$(j-1) \bigg(1 + \frac{C}{\log n} \bigg) \bigg(\frac{r}{\beta n^2 (\log n)^{1/2}} \bigg) \frac{C r^2}{n^4 \log n} \leq \frac{C (j-1) r^3}{n^6 (\log n)^{3/2}} \leq \frac{C (j-1)}{(\log n)^{6/22}}.$$
Likewise, using Proposition \ref{closeprop}, the expected value of the fourth term is at most
$C (j-1)^2 (\log n)^{-8/22}$.  

We now consider the second term.  By Proposition \ref{maincapprop}, for all $i$ we have
\begin{equation}
P \bigg( \bigg| \mbox{Cap}_{r - 2w}\big(LE(R_i) \big) - \frac{a_n}{(\log n)^{2/11}} \bigg| > \frac{1}{(\log n)^{5/22}} \bigg) \leq \frac{C}{(\log n)^{3/22 - \theta}}
\label{newcapeq}
\end{equation}
Let $V$ be the event that $|\mbox{Cap}_{r-2w}(LE(R_i)) - a_n (\log n)^{-2/11}| \leq (\log n)^{-5/22}$ for all $i \leq \ell$, and let $V_j$ be the event that this inequality holds for all $i \leq j-1$.
By (\ref{mdef}), we have $1/(m+1) \leq a_n/(\log n)^{2/11} \leq 1/m$, so $|1/m - a_n/(\log n)^{2/11}| \leq 1/m^2 \leq C/(\log n)^{4/11}$.  Therefore, on $V_j$, we have $|1/m - \mbox{Cap}_{r-2w}(LE(R_i))| \leq C/(\log n)^{5/22}$ for all $i \in S_{j-1}$.  Also, $1/(\alpha m^{3/2}) \leq C/(\log n)^{3/11}$ and, since $|E_{i,j}'| \leq C/(\log n)^{1/11}$, we have $|E_{i,j}' \mbox{Cap}_{r-2w}(LE(R_i))| \leq C/(\log n)^{3/11}$ for all $i \in S_{j-1}$ on $V_j$.  Therefore, on $V_j$, the second term on the right-hand side of (\ref{bigbound}) is at most $Cj/(\log n)^{5/22}$.

Taking expectations of both sides in (\ref{bigbound}), using that $V_j \in {\cal F}_{j-1}$, and combining these bounds, we get
\begin{equation}
P(W_j \cap V) \leq P(W_j \cap V_j) = E[P(W_j|{\cal F}_{j-1}) {\bf 1}_{V_j}] \leq \frac{C j}{(\log n)^{5/22}} + \frac{C j^2}{(\log n)^{8/22}}.
\label{finnew1}
\end{equation}
Let $K = C' (\log n)^{2/22} (\log \log n)$, where $C'$ is the constant from Lemma \ref{Pllem}.  Then, using 
(\ref{newcapeq}), (\ref{finnew1}), and Lemma \ref{Pllem}, we get
\begin{align}
P \bigg( \bigcup_{j=1}^{\ell} &W_j \bigg) \leq P \bigg( \bigcup_{j=1}^{\ell} (W_j \cap V) \bigg) + P(V^c) \nonumber \\
&\leq P(\ell > K) + \sum_{j=1}^{\lfloor K \rfloor} P(W_j \cap V) + \sum_{j=1}^{\lfloor K \rfloor} P \bigg( \bigg| \mbox{Cap}_{r-2w}(LE(R_j)) - \frac{a_n}{(\log n)^{2/11}} \bigg| > \frac{1}{(\log n)^{5/22}} \bigg) \nonumber \\
&\leq \frac{C}{\log n} + \frac{C K^2}{(\log n)^{5/22}} + \frac{C K^3}{(\log n)^{8/22}} + \frac{CK}{(\log n)^{3/22 - \theta}} \leq \frac{C (\log \log n)}{(\log n)^{1/22 - \theta}}. \nonumber
\end{align}
Since this result holds for all $\theta > 0$, the proposition follows.
\end{proof}

\subsection{Distances in the uniform spanning tree}

Recall from the introduction that $y_1, \dots, y_k$ are $k$ points picked uniformly at random from $K_m$, and $\tilde{\cal T}$ is a uniform spanning tree on $K_m$.  The distance $\tilde{d}(y_i, y_j)$ is the number of vertices on the path between $y_i$ and $y_j$ in $\tilde{\cal T}$, and the distance $\tilde{d}^*(y_i, y_j)$ is the length of the path between $y_i$ and $y_j$ in the tree $\tilde{\cal T}^*$, unless the path goes through $\rho$ in which case $\tilde{d}^*(y_i, y_j) = \infty$.  Likewise, recall that $x_1, \dots, x_k$ are points chosen uniformly at random from $\Z^4_n$.  The distance $d(x_i, x_j)$ is the number of vertices on the path between $x_i$ and $x_j$ in the uniform spanning tree ${\cal T}$ on $\Z^4_n$.  Also, $d'(x_i, x_j)$ is the length of the path from $x_i$ to $x_j$ in the tree ${\cal T}_k$, while $d^*(x_i, x_j)$ is the length of the path from $\zeta_i$ to $\zeta_j$ in the tree ${\cal T}^*$.  Both $d'(x_i, x_j)$ and $d^*(x_i, x_j)$ are set to $\infty$ when the path goes through the root.  Peres and Revelle proved the following result in section 7 of \cite{perrev}.

\begin{Lemma}
The total variation distance between the joint distribution of the $\binom{k}{2}$ distances $(d(x_i, x_j))_{1 \leq i < j \leq k}$ and the joint distribution of the distances $(d'(x_i, x_j))_{1 \leq i < j \leq k}$ is at most the probability that for some $i, j \leq k$, the path in ${\cal T}_k$ from $x_i$ to $x_j$ goes through the root $\rho$.  Likewise, the total variation distance between the distributions of $(\tilde{d}(y_i, y_j))_{1 \leq i < j \leq k}$ and $(\tilde{d}^*(y_i, y_j))_{1 \leq i < j \leq k}$ is at most the probability that for some $i, j \leq k$, the path in $\tilde{\cal T}^*$ from $y_i$ to $y_j$ goes through $\rho$.
\label{fmlem}
\end{Lemma}

The argument of Peres and Revelle relies on a form of stochastic domination of spanning forests by spanning trees, which follows from a result of Feder and Mihail \cite{fedmih}.  Note that ${\cal T}_k$ is a portion of a uniform spanning tree on $G_{n, \beta}$ (the rest of the tree can be constructed by continuing with Wilson's algorithm), and when edges in this spanning tree leading to the root $\rho$ are removed, the graph that remains is a random spanning forest on $\Z^4_n$.  See section 7 of \cite{perrev} for details.

Because the complete graph satisfies the conditions of the graphs studied by Peres and Revelle, the following Lemma follows immediately from Lemma 7.2 of \cite{perrev}.

\begin{Lemma}
Let $\epsilon > 0$.  For sufficiently large $\alpha$, there exists an $M$ such that if $m > M$ then the probability that for some $i, j \leq k$, the path in $\tilde{\cal T}^*$ from $y_i$ to $y_j$ goes through $\rho$ is less than $\epsilon$.
\label{domlem}
\end{Lemma}

The next lemma bounds the number of segments that are involved in intersections.  Recall that $J = \{j \leq \ell: I_{i,j} = 1 \mbox{ or }I_{j,i} = 1 \mbox{ for some }0 \leq i \leq \ell\}$.  Let $|J|$ denote the cardinality of $J$.

\begin{Lemma}
There exists a constant $C$ depending on $k$ such that  $$P \big( |J| > (\log n)^{1/22} \big) \leq \frac{C (\log \log n)^2}{(\log n)^{1/22}}.$$
\label{maindistcouple}  
\end{Lemma}

\begin{proof}
As in the proof of Proposition \ref{coupleprop}, let $K = C'(\log n)^{2/22} (\log \log n)$, where $C'$ is the constant from Lemma \ref{Pllem}.  By Lemma \ref{Pllem}, we have $P(\ell > K) \leq C/(\log n)$.  By taking expectations in (\ref{condIij}) and using Proposition \ref{lpssmallprop}, we get $$P(I_{i,j} = 1) \leq \frac{Cr^2}{n^4 \log n}$$ for all $1 \leq i < j$.  By (\ref{condI0j}), $P(I_{0,j} = 1) \leq Cr/(n^2 (\log n)^{1/2})$ for all $j$.  Therefore, for all $i \leq K$, $$P(I_{i,j} = 1 \mbox{ or }I_{j,i} = 1 \mbox{ for some }0 \leq j \leq K) \leq \frac{CK r^2}{n^4 \log n} + \frac{C r}{n^2 (\log n)^{1/2}} \leq \frac{C \log \log n}{(\log n)^{2/22}}.$$  It follows that $$E[|J \cap \{1, \dots, K\}|] \leq\frac{CK \log \log n}{(\log n)^{2/22}} \leq C(\log \log n)^2.$$  Using Markov's Inequality,
$$P \big( |J| > (\log n)^{1/22} \big) \leq P\big( |J \cap \{1, \dots, K\}| > (\log n)^{1/22} \big) + P(\ell > K) \leq \frac{C (\log \log n)^2}{(\log n)^{1/22}} + \frac{C}{\log n},$$
and the lemma follows.
\end{proof}

\begin{Lemma}
There exists a constant $C$ depending on $k$ and a sequence of constants $(\gamma_n)_{n=1}^{\infty}$ with $0 < \inf \gamma_n \leq \sup \gamma_n < \infty$ such that outside of an event whose probability is at most $C(\log \log n)^2 (\log n)^{-1/22}$, we have
$$\bigg| \frac{d^*(x_i, x_j)}{\sqrt{m}} - \frac{d'(x_i, x_j)}{\gamma_n n^2 (\log n)^{1/6}} \bigg| \leq \frac{C}{(\log n)^{1/22}}.$$
\label{dcompare}
\end{Lemma}

\begin{proof}
By Proposition \ref{distprop}, on $G$ we have
\begin{equation}
\bigg| d'(x_i, x_j) - \sum_{h \in S_{i,j}} N_h \bigg| \leq \sum_{h \in J} N_h.
\label{tri1}
\end{equation}
By recalling that $d^*(x_i, x_j) = |S_{i,j}| = |J \cap S_{i,j}| + |J^c \cap S_{i,j}|$ and splitting the sum into $h \in J$ and $h \in J^c$, we get
\begin{align}
\bigg| b_n n^2 (\log n)^{5/66} &d^*(x_i, x_j) - \sum_{h \in S_{i,j}} N_h \bigg| \nonumber \\
&\leq b_n n^2 (\log n)^{5/66} |J| + \sum_{h \in J} N_h + \sum_{h \in J^c \cap S_{i,j}} \big| N_h - b_n n^2 (\log n)^{5/66} \big|,
\label{tri2}
\end{align}
where $(b_n)_{n=1}^{\infty}$ is the sequence of constants defined in Corollary \ref{lengthcor} with $L = r$.  Choose $\gamma_n$ such that $\sqrt{m} b_n n^2 (\log n)^{5/66} = \gamma_n n^2 (\log n)^{1/6}$.  As $m = \lfloor a_n^{-1} (\log n)^{2/11} \rfloor$, we get $\gamma_n = b_n a_n^{-1/2} (1 - \delta_n)$, where $0 \leq \delta_n \leq C (\log n)^{-2/11}$ for some constant $C$.  The fact that $0 < \inf \gamma_n \leq \sup \gamma_n < \infty$ now follows from Corollary \ref{lengthcor} and Proposition \ref{maincapprop}.  Combining (\ref{tri1}) and (\ref{tri2}), we get
\begin{align}
\bigg| \frac{d^*(x_i, x_j)}{\sqrt{m}} - \frac{d'(x_i, x_j)}{\gamma_n n^2 (\log n)^{1/6}} \bigg| &\leq
\bigg| \frac{d^*(x_i, x_j)}{\sqrt{m}} - \frac{1}{\sqrt{m} b_n n^2 (\log n)^{5/66}} \sum_{h \in S_{i,j}} N_h \bigg| \nonumber \\
&\hspace{.3in}+ \bigg| \frac{1}{\gamma_n n^2 (\log n)^{1/6}} \sum_{h \in S_{i,j}} N_h - \frac{d'(x_i, x_j)}{\gamma_n n^2 (\log n)^{1/6}} \bigg| \nonumber \\
&\leq \frac{2}{\gamma_n n^2 (\log n)^{1/6}} \sum_{h \in J} N_h + \frac{|J|}{\sqrt{m}} \nonumber \\
&\hspace{.3in}+ \frac{1}{\gamma_n n^2 (\log n)^{1/6}} \sum_{h \in J^c \cap S_{i,j}} \big| N_h - b_n n^2 (\log n)^{5/66} \big|.
\label{3bound}
\end{align}
It remains to show that the three terms on the right-hand side of (\ref{3bound}) are small.  

Let $H$ be the event that $|LE(R_j) - b_n n^2 (\log n)^{5/66}| \leq C' n^2 (\log n)^{-15/44}$ for all $j \leq \ell$, where $C'$ is the constant from Corollary \ref{lengthcor}.  Then Corollary \ref{lengthcor} and Lemma \ref{Pllem} imply $$P(H^c) \leq P(\ell > K) + \frac{CK}{(\log n)^{4/22}} \leq \frac{C(\log \log n)}{(\log n)^{1/11}}.$$  Combining this observation with Proposition \ref{goodprop} and Lemmas \ref{Pllem} and \ref{maindistcouple}, we get $$P \big(G \cap H \cap \{|J| \leq (\log n)^{1/22} \} \cap \{\ell \leq K \} \big) \geq 1 - \frac{C (\log \log n)^2}{(\log n)^{1/22}}.$$
By Lemma \ref{labelprop}, on this event we have
$$\frac{2}{\gamma_n n^2 (\log n)^{1/6}} \sum_{h \in J} N_h \leq \frac{C |J|}{n^2 (\log n)^{1/6}} \cdot n^2 (\log n)^{5/66} \leq \frac{C}{(\log n)^{1/22}}$$ and $$\frac{|J|}{\sqrt{m}} \leq \frac{C}{(\log n)^{1/22}}.$$  Also, again using Lemma \ref{labelprop},
\begin{align}
\frac{1}{\gamma_n n^2 (\log n)^{1/6}} \sum_{h \in J^c \cap S_{i,j}} \big| N_h - b_n n^2 (\log n)^{5/66} \big| &\leq \frac{C \ell}{n^2 (\log n)^{1/6}} \bigg( 2w + \frac{C' n^2}{(\log n)^{15/44}} \bigg) \nonumber \\
&\leq \frac{C (\log \log n)}{(\log n)^{1/6 - 1/11 - \theta}} = \frac{C (\log \log n)}{(\log n)^{5/66 - \theta}}. \nonumber
\end{align}
The result follows from (\ref{3bound}) and these bounds.
\end{proof}

\begin{proof}[Proof of Theorem \ref{mainth}]
By Aldous' result (\ref{aldCG}), the joint distribution of $\tilde{d}(y_i, y_j)/\sqrt{m}$ for $1 \leq i < j \leq k$ converges to $\mu_k$ as $m \rightarrow \infty$.  Therefore, given $\epsilon > 0$, it suffices to show that for sufficiently large $n$, the total variation distance between the joint distribution of these distances and the joint distribution of the distances $d(x_i, x_j)/(\gamma_n n^2 (\log n)^{1/6})$ is less than $\epsilon$.

By (\ref{alphadef}) Lemma \ref{domlem}, we can choose $\alpha$ and $\beta$ large enough that for sufficiently large $m$, the probability that for some $i, j \leq k$ the path in $\tilde{\cal T}^*$ from $y_i$ to $y_j$ goes through $\rho$ is less than $\epsilon$.  Therefore, by Lemma \ref{fmlem}, the total variation distance between the distributions of $(\tilde{d}(y_i, y_j))_{1 \leq i < j \leq k}$ and $(\tilde{d}^*(y_i, y_j))_{1 \leq i < j \leq k}$ is less than $\epsilon$ for sufficiently large $m$, and thus for sufficiently large $n$.

Suppose ${\cal T}^* = \tilde{\cal T}^*$.  If for all $i,j \leq k$, the path in $\tilde{\cal T}^*$ from $y_i$ to $y_j$ does not go through $\rho$, then for all $i,j \leq k$, the path in ${\cal T}^*$ from $\zeta_i$ to $\zeta_j$ does not go through $\rho$.  However, on $G$ this happens if and only if for all $i, j \leq k$, the path in ${\cal T}_k$ from $x_i$ to $x_j$ does not go through $\rho$.  In view of Proposition \ref{goodprop} and Remark \ref{Trem}, it follows that for sufficiently large $n$, the total variation distance between the distributions of $(d(x_i, x_j))_{1 \leq i < j \leq k}$ and $(d'(x_i, x_j))_{1 \leq i < j \leq k}$ is less than $\epsilon$.

For $i, j \leq k$, we have $d^*(x_i, x_j) = \tilde{d}^*(y_i, y_j)$ on the event ${\cal T}^* = \tilde{\cal T}^*$.  Thus, for sufficiently large $n$, the total variation distance between the distributions of $(\tilde{d}(y_i, y_j))_{1 \leq i < j \leq k}$ and $(d^*(x_i, x_j))_{1 \leq i < j \leq k}$ is less than $\epsilon$.  The result follows from these observations and Lemma \ref{dcompare}.
\end{proof}

\bigskip
\smallskip
\noindent {\bf \Large Acknowledgments}

\bigskip
\noindent The author thanks Yuval Peres, Jim Pitman, and David Revelle for helpful discussions, and Gady Kozma for comments on a previous version.  He also thanks a referee for suggestions which led to simplifications of several proofs in the paper.

\end{document}